\input amstex
\documentstyle{amsppt}

\define\Cee{{\Bbb C}}

\define\Pee{\Bbb P}

\define\Zee{\Bbb Z}

\define\Res#1{\operatorname{Res}_{#1 = 0}}
\define\hil#1{{X^{[#1]}}}

\define\w{\tilde}

\def\sumlim#1#2{\underset {#1} \to {\overset {#2} \to \sum}}

\def\bigopluslim#1#2{\underset {#1} \to {\overset {#2} \to \bigoplus}}
\def\exact#1#2#3{0\rightarrow{#1}\rightarrow{#2}\rightarrow{#3}\rightarrow0}

\define\End{\operatorname{End}}

\define\Id{\operatorname{Id}}

\define\Supp{\operatorname{Supp}}
\define\ad{\operatorname{ad}}
\define\ch{\operatorname{ch}}
\define\td{\operatorname{td}}

\define\proof{\demo{Proof}}
\def\stop {\nobreak$\quad$\lower 1pt\vbox{
    \hrule
    \hbox to 7pt{\vrule height 7pt\hfil\vrule height 7pt}
      \hrule}\ifmmode\relax\else\par\medbreak \fi}
\define\endproof{\stop\enddemo}
\define\endstatement{\rm }
\define\theorem#1{\medskip\noindent {\bf {Theorem #1}.} \it}
\define\lemma#1{\medskip\noindent {\bf {Lemma #1}.} \it}
\define\proposition#1{\medskip\noindent {\bf {Proposition #1}.} \it}
\define\corollary#1{\medskip\noindent {\bf {Corollary #1}.} \it}
\define\claim#1{\medskip\noindent {\bf {Claim #1}.} \it}

\define\section#1{\bigskip\noindent {\bf #1}}
\define\ssection#1{\medskip\noindent{\bf #1}}

\def\today{\ifcase \month \or January\or February\or
  March\or April\or May\or June\or July\or August\or
  September\or October\or November\or December\fi\
  \number \day, \number \year}

\loadbold
\hoffset=0.3 in
\voffset=0.6 in

\leftheadtext{Wei-ping Li, Zhenbo Qin and Weiqiang Wang}
\rightheadtext{Vertex algebras and Hilbert schemes}

\topmatter
\title
Vertex algebras and the cohomology ring structure of
Hilbert schemes of points on surfaces
\endtitle
\author {Wei-ping Li$^1$, Zhenbo Qin$^2$ and Weiqiang Wang$^3$}
\endauthor
\address
Department of Mathematics, HKUST, Clear Water Bay, Kowloon, Hong Kong
\endaddress
\email mawpli\@uxmail.ust.hk
\endemail
\address
Department of Mathematics, University of Missouri, Columbia,
MO 65211, USA
\endaddress
\email
zq\@math.missouri.edu
\endemail
\address
Department of Mathematics, North Carolina State University,
Raleigh, NC 27695, USA
\email
ww9c\@weyl.math.virginia.edu
\endemail
\endaddress

\address
Current Address for W. Wang: Department of Mathematics, University of
Virginia, Charlottesville, VA 22904
\endaddress 
 \email ww9c\@virginia.edu
\endemail

\thanks
${}^1$Partially supported by the grant HKUST6170/99P
\endthanks
\thanks
${}^2$Partially supported by an NSF grant and
an Alfred P. Sloan Research Fellowship
\endthanks
\thanks
${}^3$Partially supported by an NSF grant and an FR\&PD grant at NCSU
\endthanks
\abstract Using vertex algebra techniques, we determine a set of
generators for the cohomology ring of the Hilbert schemes of
points on an arbitrary smooth projective surface over the field of
complex numbers.
\endabstract
\endtopmatter

\NoBlackBoxes
\TagsOnRight
\document
\section{1. Introduction}

The Hilbert scheme $\hil{n}$ of points on a smooth projective
surface $X$ is a desingularization of the $n$-th symmetric product
of $X$ (see \cite{Fog}). An element $\xi$ in $\hil{n}$ is a
length-$n$ $0$-dimensional closed subscheme of $X$. Recently,
there are two surprising discoveries, mainly due to the work of
G\"ottsche \cite{Go1}, Nakajima \cite{Na1, Na2} and Grojnowski
\cite{Gro}, that the sum of the cohomology groups $\Bbb H_n =
H^*(X^{[n]})$ with $\Bbb Q$-coefficients of the Hilbert schemes
$\hil{n}$ for $n \ge 0$ have relationships with modular forms on
the one hand and with representations of infinite dimensional
Heisenberg algebras on the other hand (see aslo the work of Vafa
and Witten \cite{V-W} for connections with string theory). These
results have been used by Lehn \cite{Leh} to investigate the
relation between the Heisenberg algebra structure and the cup
product structure of $\Bbb H_n$. In particular, Lehn constructed
the Virasoro algebra in a geometric fashion and studied certain
tautological sheaves over $\hil{n}$. Some other recent work on
Hilbert schemes includes \cite{dCM, EGL, Go2, Hai, K-T, LQZ, Wan}.

In this paper, by using vertex algebra techniques
(see \cite{Bor, FLM, Kac}) and generalizing the work of
Nakajima, Grojnowski and Lehn \cite{Na1, Gro, Na2, Leh},
we study the cohomology ring  structure of the
Hilbert schemes
$\hil{n}$. We determine the ring generators of $H^*(X^{[n]})$  for an
arbitrary smooth projective surface
$X$ over the field of complex numbers. In particular,
we recover the result of Ellingsrud and Str\o mme \cite{ES2}
for $X = \Pee^2$. More precisely, we find a set of
$(n \cdot \dim H^*(X))$ generators for the cohomology ring $\Bbb H_n$,
and interpret the relations among these generators in terms of
certain operators in $\End (\Bbb H)$ where
$\Bbb H = \bigoplus_{n \ge 0} \Bbb H_n$. Our results also clearly
indicate that there are deep interplays between the geometry of
Hilbert schemes and vertex algebra structures which go beyond
the Heisenberg and Virasoro algebras.

To state our result,
we establish some notations and refer the details to Definition 5.1.
Let $\Cal Z_n$ be the universal codimension-$2$ subscheme
of $\hil{n} \times X$, and $p_1$ and $p_2$ be the projections of
$\hil{n} \times X$ to $\hil{n}$ and $X$ respectively.
For $\gamma \in H^s(X)$ and $n \ge 0$, let $G_i(\gamma, n)$
be the $H^{s+2i}(\hil{n})$-component of
$$G(\gamma, n) \overset \text{def} \to =
p_{1*}(\ch(\Cal O_{\Cal Z_n}) \cdot p_2^*\td(X)
\cdot p_2^*\gamma) \in \Bbb H_n  \tag 1.1$$
(we refer to the Conventions at the end of this section
for the conventions used in the paper).
For $\gamma \in H^*(X)$ and $i \in \Zee$,
define an operator $\frak G_i(\gamma) \in \End({\Bbb H})$
which acts on the component $\Bbb H_n = H^*(\hil{n})$
by the cup product by the class $G_i(\gamma, n)$.

\theorem{1.2} Let $X$ be a smooth projective surface over the field of
complex numbers. For $n \ge 1$, the cohomology ring
$\Bbb H_n = H^*(\hil{n})$ is generated by
$$G_{i}(\gamma, n) = \frak G_{i}(\gamma)(1_{\hil{n}})$$
where $0 \le i < n$ and $\gamma$ runs over a linear basis of $H^*(X)$.
Moreover, the relations among these generators
are precisely the relations among the restrictions
$\frak G_i(\gamma)|_{\Bbb H_n}$ of the corresponding operators
$\frak G_i(\gamma)$ to ${\Bbb H_n}$.
\endstatement

The above Theorem and its proof are inspired mainly by two sources.
The first one is Lehn's approach of determining the cohomology ring
structure of $(\Cee^2)^{[n]}$ by using vertex operator techniques
(the cohomology ring structure of $(\Cee^2)^{[n]}$
has been first obtained by Ellingsrud and Str\o mme \cite{ES2}).
Lehn's approach is very instrumental and valuable to us.
Our first result here is that although it is difficult to describe
completely the operators $\frak G_i(\gamma)$ as vertex operators or
differential operators, we are able to determine the leading terms
of the operators $\frak G_i(\gamma)$ as the degree-$0$ components of
some explicit vertex operators (see the paragraph preceding
Theorem 4.12 for the definition of the leading term).
These vertex operators are natural
generalization of the Virasoro operators $\frak L_n(\alpha)$
considered by Lehn \cite{Leh}. Such descriptions of the leading
terms allow us to use induction to derive our Theorem above.
As a byproduct, we also show that the commutator between
the operator $\frak G_i(\gamma)$ and the Heisenberg generator
$\frak q_n(\alpha)$ depends only on the cup product $\gamma \alpha$,
which we refer to as the {\it transfer property}. We
remark that such a transfer property seems to be
a general phenomenon among this type of commutation relations.

The second one is the work of Ellingsrud and Str\o mme \cite{ES1, ES2}
(see also the work of Beauville, Fantechi, G\"ottsche,
Yoshioka and Markman \cite{Bea, F-G, Yos, Mar}
on the cohomology ring structures of other moduli spaces of sheaves).
In \cite{ES2}, Ellingsrud and Str\o mme proved that for $X = \Pee^2$,
the cohomology ring $\Bbb H_n = H^*(\hil{n})$ is generated by
the Chern classes of the tautological rank-$n$ bundles
$$p_{1*}(\Cal O_{\Cal Z_n} \otimes p_2^*\Cal O_{\Pee^2}(-j))  \tag 1.3$$
with $j = 1, 2, 3$. Equivalently, this says that
$\Bbb H_n$ is generated by the $H^{2i}(\hil{n})$-components
with $0 \le i \le n$ of
the Chern characters of the three bundles in (1.3).
By the Grothendieck-Riemann-Roch Theorem \cite{Har}, we have
$$\ch(p_{1*}(\Cal O_{\Cal Z_n} \otimes p_2^*\Cal O_{\Pee^2}(-j)))
= p_{1*}(\ch(\Cal O_{\Cal Z_n} \otimes p_2^*\Cal O_{\Pee^2}(-j))
\cdot p_2^*\td(X))$$
for $j = 1, 2, 3$. Note that
$\ch(\Cal O_{\Cal Z_n} \otimes p_2^*\Cal O_{\Pee^2}(-j))
= \ch(\Cal O_{\Cal Z_n}) \cdot p_2^*\ch(\Cal O_{\Pee^2}(-j))
= \ch(\Cal O_{\Cal Z_n}) \cdot p_2^*([X] -j[\ell]+j^2[x]/2)$
where $\ell$ (respectively, $x$) stands for a line
(respectively, a point) in $X = \Pee^2$. So the result of
Ellingsrud and Str\o mme says that the cohomology ring $\Bbb H_n$
is generated by the $H^{2i}(\hil{n})$-components of
$$p_{1*}(\ch(\Cal O_{\Cal Z_n}) \cdot p_2^*\td(X) \cdot p_2^*\gamma)
\in \Bbb H_n$$
where $0 \le i \le n$ and $\gamma = [X], [\ell], [x]$.
Indeed, this is our motivation for (1.1).

We may also put our present work in a different perspective.
>From the pioneering work \cite{Na1, Gro, Leh}, it is clear that
there are deep connections between the geometry of Hilbert schemes
and vertex operators. However, the structure of a vertex algebra
\cite{Bor, FLM, Kac} is far more richer than the appearance of the
Heisenberg and Virasoro vertex operators which are of conformal
weight one and two respectively. Thus a natural question here,
which is easy to post but difficult to answer, is to understand
the full symmetry of vertex algebras in terms of the geometry of
Hilbert schemes. Our present work provides a strong evidence that
the vertex operators of higher conformal weights afford
nice geometric interpretations. In a work in progress,
we will further clarify precise connections between
vertex algebras and the geometry of Hilbert schemes.

The paper is organized as follows. In section 2,
we recall constructions and results of
Nakajima, Grojnowski and Lehn.
In section 3, we study the relation between
Lehn's boundary operator and the Heisenberg generators,
and introduce the transfer property.
In section 4, we prove that the leading terms for
certain linear operators of geometric significance are
the degree-$0$ components of some vertex operators. Finally,
we prove our main results in section 5.

\medskip\noindent
{\bf Conventions:} Throughout the paper, all cohomology groups
are in $\Bbb Q$-coefficients. The cup product between
two cohomology classes $\alpha$ and $\beta$
is denoted by $\alpha \cdot \beta$ or simply by $\alpha \beta$.
For a continuous map $p: Y_1 \to Y_2$ between two smooth compact
manifolds and for $\alpha_1 \in H^*(Y_1)$, the push-forward
$p_*(\alpha_1)$ is defined by
$$p_*(\alpha_1) = \text{PD}^{-1}p_{*}(\text{PD}(\alpha_1))$$
where $\text{PD}$ stands for the Poincar\'e duality.
We make no distinction between an algebraic cycle and
its corresponding cohomology class so that intersections among
algebraic cycles correspond to cup products among
the corresponding cohomology classes. For instance, for two algebraic
cycles $[a]$ and $[b]$ on a smooth projective variety $Y$,
it is understood that $[a] \cdot [b] \in H^*(Y)$.

\bigskip\noindent
{\bf Acknowledgments:} After communicating this paper to M. Lehn,
we were informed by him that he and C. Sorger had a preprint on
the cohomology ring structure of $\hil{n}$ for $X = \Cee^2$,
and that their methods could apply to the case when $K_X = 0$ by
 combining with the results in the paper (see \cite{LS1, LS2}).
The authors thank J. Li, Y. Ruan and Y. Zhu for stimulating discussions,
and M. Lehn for helpful comments. Also, the first author thanks
the University of Missouri at Columbia for its hospitality
during his visit in August 2000, and the second and third authors thank
the Hong Kong University of Science and Technology for its hospitality
and financial support during their visits in July 2000.

\section{2. Results of Nakajima, Grojnowski and Lehn}

In this section, we shall fix some notations,
and recall some results of Nakajima, Grojnowski and Lehn
\cite{Na1, Gro, Leh}.
For convenience, we also review certain basic facts for
the Hilbert scheme of points in a smooth projective surface.

Let $X$ be a smooth projective surface over $\Cee$,
and $\hil{n}$ be the Hilbert scheme of points in $X$.
An element in the Hilbert scheme $\hil{n}$ is represented
by a length-$n$ $0$-dimensional
closed subscheme $\xi$ of $X$, which sometimes is called
a length-$n$ $0$-cycle. For $\xi \in \hil{n}$, let $I_{\xi}$ and
$\Cal O_\xi$ be the corresponding sheaf of ideals and
structure sheaf respectively.
For a point $x \in X$, let $\xi_x$ be the component of $\xi$
supported at $x$ and $I_{\xi, x} \subset \Cal O_{X, x}$ be
the stalk of $I_{\xi}$ at $x$.
It is known from \cite{Fog} that $\hil{n}$ is smooth.
In $\hil{n}\times X$, we have the universal codimension-$2$ subscheme:
$$\Cal Z_n=\{(\xi, x) \subset \hil{n}\times X \, | \,
x\in \Supp{(\xi)}\}\subset \hil{n}\times X.$$
In general, for $n > m > 0$, we have the closed subscheme
$$X^{[n, m]} = \{(\xi, \eta) \in \hil{n}\times \hil{m} \, | \,
\xi \supset \eta \}$$
of $\hil{n}\times \hil{m}$. Let $\psi$ and $\phi$ be the natural maps
from $X^{[n, m]}$ to $\hil{n}$ and $\hil{m}$ respectively.
Then, $X^{[n, m]}$ parametrizes the two flat families
$\psi_{X}^{-1}(\Cal Z_n) \supset \phi_{X}^{-1}(\Cal Z_m)$
over $X^{[n, m]} \times X$ where
$\psi_{X} = \psi \times \Id_X: X^{[n, m]} \times X \to X^{[n]} \times X$
and $\phi_{X} = \phi \times \Id_X: X^{[n, m]} \times X
\to X^{[m]} \times X$. In addition, there exists an exact sequence
$$\exact{\Cal I_{n, m}}{\psi_{X}^*\Cal O_{\Cal Z_n}}
{\phi_{X}^*\Cal O_{\Cal Z_m}}. \tag 2.1$$
When $m = n-1$, there is a morphism $\rho: \hil{n,n-1} \to X$
which maps a point $(\xi, \eta) \in X^{[n, n-1]}$
to the support of $(I_\eta/I_\xi)$. So we have the standard diagram
$$\matrix
X&\overset \rho \to \leftarrow&X^{[n, n-1]}&\overset \psi \to \to&\hil{n}\\
 &                        &\downarrow{\phi}&                     &       \\
 &                              &\hil{n-1}.&                     &       \\
\endmatrix \tag 2.2$$
It is known from \cite{Che,Tik, ES3} that
$\hil{n,n-1}$ is irreducible, smooth and of dimension $2n$.
In fact, $\hil{n,n-1}$ is isomorphic to the blowup of
$\hil{n-1}\times X$ along $\Cal Z_{n-1}$. Let $E_n$ be
the exceptional divisor in $\hil{n,n-1}$. Then, we have
$$E_n = \{(\xi, \eta) \in \hil{n, n-1} \, | \,
\Supp(\xi) = \Supp(\eta) \}.$$
Moreover, (2.1) can be simplified to the exact sequence
(see p.193 in \cite{Leh}):
$$\exact{\rho_X^*\Cal O_{\Delta_X} \otimes
p_1^*\Cal O_{\hil{n,n-1}}(-E_n)}{\psi_{X}^*\Cal O_{\Cal Z_n}}
{\phi_{X}^*\Cal O_{\Cal Z_{n-1}}}   \tag 2.3$$
where $\Delta_X$ is the diagonal in $X \times X$, and $p_1$ is
the projection of $X^{[n, n-1]} \times X$ to $X^{[n, n-1]}$.
Finally, we let $X^n=\underbrace{X\times\cdots\times X}_{n}$
be the $n$-th Cartesian product, and
$$X^{[n_1], \ldots, [n_k]} = \hil{n_1} \times \cdots \times \hil{n_k}.
\tag 2.4$$

We formulate below various notations and definitions which
will be used later.

\medskip\noindent
{\bf Definition 2.5.} (i) Let $\Bbb H =
\bigopluslim{n, i \ge 0}{} \Bbb H^{n,i}$
denote the double graded vector space with components
$\Bbb H^{n,i} \overset \text{def} \to = H^i(\hil{n})$, and
$\Bbb H_n \overset \text{def} \to = H^*(\hil{n})
\overset \text{def} \to = \bigopluslim{i=0}{4n} H^i(\hil{n})$.
The element $1$ in $H^0(\hil{0}) = \Bbb Q$
is called the {\it vacuum vector} and denoted by $|0\rangle$;

(ii) A linear operator $\frak f \in \End(\Bbb H)$ is
{\it homogeneous of bi-degree} $(\ell, m)$ if
$$\frak f(\Bbb H^{n,i}) \subset \Bbb H^{n+\ell,i+m}.  \tag 2.6$$
Furthermore, $\frak f \in \End(\Bbb H)$ is
{\it even} (respectively, {\it odd}) if $m$ is even (respectively, odd).

(iii) For two homogeneous linear operators
$\frak f$ and $\frak g \in \End(\Bbb H)$ of bi-degrees
$(\ell, m)$ and $(\ell_1, m_1)$ respectively,
define the {\it Lie superalgebra bracket} $[\frak f, \frak g]$ by
$$[\frak f, \frak g] = \frak f \frak g -
(-1)^{m m_1} \frak g \frak f.   \tag 2.7$$

A non-degenerate super-symmetric bilinear form $(, )$ on $\Bbb H$
is induced from the standard one on $\Bbb H_n=H^*(\hil{n})$.
For a homogeneous linear operator $\frak f \in \End(\Bbb H)$ of
bi-degree $(\ell, m)$, we can define its {\it adjoint}
$\frak f^\dagger \in \End(\Bbb H)$ by
$$(\frak f(\alpha), \beta) = (-1)^{m \cdot |\alpha|}
\cdot (\alpha, \frak f^\dagger(\beta))$$
where $|\alpha| = s$ for $\alpha \in H^s(X^{[n]})$. Note that the
bi-degree of $\frak f^\dagger$ is $(-\ell, m - 4 \ell)$. Also,
$$(\frak f \frak g)^\dagger = (-1)^{m m_1} \cdot
\frak g^\dagger \frak f^\dagger \qquad \text{and} \qquad
[\frak f, \frak g]^\dagger = -[\frak f^\dagger, \frak g^\dagger]
\tag 2.8$$
where $\frak g \in \End(\Bbb H)$ is another homogeneous linear operator
of bi-degree $(\ell_1, m_1)$.

Next, we collect from \cite{Na1, Leh} the definitions of
the closed subset $Q^{[n+\ell,n]}$ in
$\hil{n+\ell} \times X \times \hil{n}$,
the Heisenberg generator $\frak q_n$,
the Virasoro generator $\frak L_n$,
the boundary operator $\frak d$, and the derivative $\frak f'$
of a linear operator $\frak f \in \End(\Bbb H)$.

\medskip\noindent
{\bf Definition 2.9.} (see \cite{Na1, Leh})
(i) For $n \ge 0$, define $Q^{[n,n]} = \emptyset$.
For $n \ge 0$ and $\ell > 0$, define $Q^{[n+\ell,n]}
\subset \hil{n+\ell} \times X \times \hil{n}$ to be the closed subset
$$\{ (\xi, x, \eta) \in \hil{n+\ell} \times X \times \hil{n} \, | \,
\xi \supset \eta \text{ and } \Supp(I_\eta/I_\xi) = \{ x \} \};
\tag 2.10$$

(ii) For $n \in \Zee$, define linear maps
$\frak q_n: H^*(X) \to \End(\Bbb H)$ as follows.
When $n \ge 0$, the linear operator $\frak q_n(\alpha)
\in \End(\Bbb H)$ with $\alpha \in H^*(X)$ is defined by
$$\frak q_n(\alpha)(a) = \w p_{1*}([Q^{[m+n,m]}]
\cdot \w \rho^*\alpha \cdot \w p_2^*a)  \tag 2.11$$
for all $a \in \Bbb H_m = H^*(\hil{m})$,
where $[Q^{[m+n,m]}]$ is (the cohomology class corresponding to)
the algebraic cycle associated to $Q^{[m+n,m]}$,
and $\w p_1, \w \rho, \w p_2$ are the projections of
$\hil{m+n} \times X \times \hil{m}$ to
$\hil{m+n}, X, \hil{m}$ respectively. When $n < 0$,
define the operator $\frak q_n(\alpha) \in \End(\Bbb H)$
with $\alpha \in H^*(X)$ by
$$\frak q_n(\alpha) = (-1)^n \cdot \frak q_{-n}(\alpha)^\dagger;
\tag 2.12$$

(iii) For $n \in \Zee$, define linear maps
$\frak L_n: H^*(X) \to \End(\Bbb H)$ by putting
$$\frak L_n =
\cases
  {1 \over 2} \cdot \sumlim{m \in \Zee}{} \frak q_m \frak q_{n-m}
     \tau_{2*}, &\text{if $n \ne 0$}  \\
  \sumlim{m > 0}{} \frak q_m \frak q_{-m} \tau_{2*}, &\text{if $n = 0$} \\
\endcases
\tag 2.13$$
where $\tau_{2*}: H^*(X) \to H^*(X^2)$ is the linear map induced by
the diagonal embedding $\tau_2: X \to X^2$,
and the operator $\frak q_m \frak q_{\ell} \tau_{2*}(\alpha)$ stands for
$$\sum_j \frak q_m(\alpha_{j,1}) \frak q_{\ell}(\alpha_{j,2})$$
when $\tau_{k*}\alpha = \sum_j \alpha_{j,1} \otimes
\alpha_{j, 2}$ via the K\"unneth decomposition of $H^*(X^2)$;

(iv) Define the linear operator $\frak d \in \End(\Bbb H)$ by
$$\frak d = \bigopluslim{n}{} c_1(p_{1*}\Cal O_{\Cal Z_n})
= \bigopluslim{n}{} (-[\partial \hil{n}]/2)   \tag 2.14$$
where $p_1$ is the projection of $\hil{n} \times X$ to $\hil{n}$,
$\partial \hil{n}$ is the boundary of $\hil{n}$ consisting of
all $\xi \in \hil{n}$ with $|\Supp(\xi)| < n$,
and the first Chern class $c_1(p_{1*}\Cal O_{\Cal Z_n})$ of
the rank-$n$ bundle $p_{1*}\Cal O_{\Cal Z_n}$
acts on $\Bbb H_n = H^*(\hil{n})$ by the cup product.

(v) For a linear operator $\frak f \in \End(\Bbb H)$,
define its {\it derivative} $\frak f'$ by
$$\frak f' \overset \text{def} \to = \ad(\frak d) \frak f
\overset \text{def} \to = [\frak d, \frak f].  \tag 2.15$$
The higher derivative $\frak f^{(k)}$ of $\frak f$ is defined
inductively by $\frak f^{(k)} = [\frak d, \frak f^{(k-1)}]$.

We remark that the definition of the Virasoro generator $\frak L_n$
will be generalized in Definition 4.3 (ii) below. Also,
$\frak q_n(\alpha)$, $\frak L_n(\alpha)$, and $\frak d$
are homogeneous of bi-degrees $(n, 2n-2+|\alpha|)$,
$(n, 2n+|\alpha|)$, and $(0, 2)$ respectively.

Finally, we recall from \cite{Na1, Leh} the formulas for
the derivative $\frak q_n'(\alpha)$ as well as the commutation relations
among the Heisenberg generators $\frak q_n(\alpha)$ and
the Virasoro generator $\frak L_m(\beta)$.
These formulas will be used frequently in the sequel.

\theorem{2.16}
Let $K_X$ and $c_2(X)$ be the canonical divisor and
the second Chern class of $X$ respectively. Let $n,m \in \Zee$
and $\alpha, \beta \in H^*(X)$. Then,
\medskip
{\rm (i)} $[\frak q_n(\alpha), \frak q_m(\beta)]
= n \cdot \delta_{n+m} \cdot \int_X(\alpha \beta) \cdot \Id_{\Bbb H}$;
\medskip
{\rm (ii)} $[\frak L_n(\alpha), \frak q_m(\beta)]
= -m \cdot \frak q_{n+m}(\alpha \beta)$;
\medskip
{\rm (iii)} $[\frak L_n(\alpha), \frak L_m(\beta)]
= (n-m) \cdot \frak L_{n+m}(\alpha \beta) - {n^3-n \over 12} \cdot
\delta_{n+m} \cdot \int_X(c_2(X) \alpha \beta) \cdot \Id_{\Bbb H}$;
\medskip
{\rm (iv)} $\frak q_n'(\alpha) = n \cdot \frak L_n(\alpha)
+ {n(|n|-1) \over 2} \frak q_n(K_X \alpha)$;
\medskip
{\rm (v)} $[\frak q_n'(\alpha), \frak q_m(\beta)]
= -nm \cdot \left \{ \frak q_{n+m}(\alpha \beta) +
{|n|-1 \over 2} \cdot \delta_{n+m} \cdot
\int_X(K_X \alpha \beta) \cdot \Id_{\Bbb H} \right \}$.
\endstatement

We notice that Theorem 2.16 (i) was proved by Nakajima \cite{Na1}
subject to some universal nonzero constant,
which was determined subsequently in [ES3].
The other four formulas in Theorem 2.16 were obtained by Lehn \cite{Leh}.
Moreover, as observed by Nakajima and Grojnowski in \cite{Na1, Gro},
$\Bbb H$ is an irreducible representation of
the Heisenberg algebra generated by
the $\frak q_i(\alpha)$'s with $|0\rangle \in H^0(X^{[0]})$
being the highest weight vector. In particular,
a linear basis of $\Bbb H$ is given by
$$\frak q_{i_1}(\alpha_1) \frak q_{i_2}(\alpha_2)
\cdots \frak q_{i_k}(\alpha_k) |0\rangle,$$
where $k \ge 0$, $i_1 \geq i_2 \geq \cdots \geq i_k >0$, and each
of the cohomology classes $\alpha_1, \alpha_2, \ldots, \alpha_k$
runs over a fixed linear basis of
$H^*(X) = \bigopluslim{i=0}{4} H^i(X)$.

\section{3. The higher derivatives of $\frak q_n(\alpha)$}

In this section, we study the higher derivatives $\frak q_n^{(k)}(\alpha)$
of the Heisenberg generator $\frak q_n(\alpha)$ by computing
their commutators with other Heisenberg generators $\frak q_m(\beta)$.
In addition, for two series of operators $\frak A(\alpha), \frak B(\beta)
\in \End(\Bbb H)$ depending linearly on the cohomology classes
$\alpha, \beta \in H^*(X)$,
we introduce the transfer property for the commutators
$[\frak A(\alpha), \frak B(\beta)]$, i.e., the property that
$$[\frak A(\alpha), \frak B(\beta)] = [\frak A(1_X), \frak B(\alpha \beta)]
= [\frak A(\alpha \beta), \frak B(1_X)].$$
Then, we prove that the commutators $[\frak q_n^{(k)}(\alpha),
\frak q_m(\beta)]$ satisfy the transfer property.

\lemma{3.1} Let $k \ge 0$, $n_0, \ldots, n_k \in \Zee$,
and $\alpha_0, \ldots, \alpha_k \in H^*(X)$. Then,
$$\align
&\qquad [ \ldots [\frak q^{(k)}_{n_0}(\alpha_0), \frak q_{n_1}(\alpha_1)],
   \ldots ], \frak q_{n_k}(\alpha_k)] \\
&=a \cdot \frak q_{n_0 + \ldots + n_k}(\alpha_0 \ldots \alpha_k)
   +b \cdot \int_X(K_X \alpha_0 \ldots \alpha_k) \cdot \Id_{\Bbb H}
   \tag 3.2  \\
\endalign$$
where $a$ and $b$ are constants depending only on $k, n_0, \ldots, n_k$.
\endstatement
\proof
For simplicity, denote the left-hand-side of (3.2) by
$\Cal C(k; n_0, \alpha_0; \ldots ; n_k, \alpha_k)$.
Note that (3.2) is trivially true for $k=0$. By Theorem 2.16 (v),
(3.2) is true for $k=1$. In the following,
by assuming that (3.2) holds for some positive integer $k$,
we shall prove that (3.2) also holds for $(k+1)$.
Recall the Jacobi identity
$$[[\frak f_1, \frak f_2], \frak f_3]
=(-1)^{m_2 m_3}[[\frak f_1, \frak f_3], \frak f_2]
+[\frak f_1,[\frak f_2, \frak f_3]]     \tag 3.3$$
if the bi-degree of $\frak f_i \in \End(\Bbb H)$ is $(\ell_i, m_i)$.
We have by the Jacobi identity that
$$\align
&\qquad \Cal C(k+1; n_0, \alpha_0; \ldots ; n_{k+1}, \alpha_{k+1}) \\
&=[ \ldots [[\frak q^{(k+1)}_{n_0}(\alpha_0), \frak q_{n_1}(\alpha_1)],
  \frak q_{n_2}(\alpha_2)], \ldots ], \frak q_{n_{k+1}}(\alpha_{k+1})] \\
&=[ \ldots [\{[\frak q^{(k)}_{n_0}(\alpha_0), \frak q_{n_1}(\alpha_1)]'
  -[\frak q^{(k)}_{n_0}(\alpha_0), \frak q_{n_1}'(\alpha_1)]\},
  \frak q_{n_2}(\alpha_2)], \ldots ], \frak q_{n_{k+1}}(\alpha_{k+1})] \\
&=[ \ldots [[\frak q^{(k)}_{n_0}(\alpha_0), \frak q_{n_1}(\alpha_1)]',
  \frak q_{n_2}(\alpha_2)], \ldots ], \frak q_{n_{k+1}}(\alpha_{k+1})] \\
&\quad - [ \ldots [\frak q^{(k)}_{n_0}(\alpha_0), \frak q_{n_1}'(\alpha_1)],
  \frak q_{n_2}(\alpha_2)], \ldots ], \frak q_{n_{k+1}}(\alpha_{k+1})].
  \tag 3.4   \\
\endalign$$
Repeating the above process, we conclude that
$$\align
&\qquad \Cal C(k+1; n_0, \alpha_0; \ldots ; n_{k+1}, \alpha_{k+1}) \\
&= [\Cal C(k; n_0, \alpha_0; \ldots ; n_k, \alpha_k),
  \frak q_{n_{k+1}}(\alpha_{k+1})]' - \sum_{i=1}^{k+1}
  [ \ldots [\frak q^{(k)}_{n_0}(\alpha_0), \frak q_{n_1}(\alpha_1)],
  \ldots      \\
&\qquad
  \ldots ], q_{n_{i-1}}(\alpha_{i-1})],
  q_{n_i}'(\alpha_i)], q_{n_{i+1}}(\alpha_{i+1})], \ldots ],
  \frak q_{n_{k+1}}(\alpha_{k+1})].  \tag 3.5   \\
\endalign$$
By the induction hypothesis on $k$ and Theorem 2.16 (i), we get
$$[\Cal C(k; n_0, \alpha_0; \ldots ; n_k, \alpha_k),
  \frak q_{n_{k+1}}(\alpha_{k+1})]'= 0.   \tag 3.6 $$
Using the Jacobi identity (3.3) again, we see that
$$\align
 &[ \ldots [\frak q^{(k)}_{n_0}(\alpha_0), \frak q_{n_1}(\alpha_1)],
  \ldots ], q_{n_{i-1}}(\alpha_{i-1})], q_{n_i}'(\alpha_i)],
  q_{n_{i+1}}(\alpha_{i+1})], \ldots ], \frak q_{n_{k+1}}(\alpha_{k+1})] \\
=& (-1)^{\sumlim{\ell =i+1}{k+1} |\alpha_\ell||\alpha_i|} \
  [\Cal C(k; n_0, \alpha_0; \ldots ; n_{i-1}, \alpha_{i-1};
  n_{i+1}, \alpha_{i+1}; \ldots ; n_{k+1}, \alpha_{k+1}),
  q_{n_i}'(\alpha_i)]  \\
 &+ [ \ldots [\frak q^{(k)}_{n_0}(\alpha_0), \frak q_{n_1}(\alpha_1)],
  \ldots ], q_{n_{i-1}}(\alpha_{i-1})],
  [q_{n_i}'(\alpha_i), q_{n_{i+1}}(\alpha_{i+1})]], \\
 &\qquad\qquad\quad q_{n_{i+2}}(\alpha_{i+2})],
  \ldots ], \frak q_{n_{k+1}}(\alpha_{k+1})]  \\
 &+ \sum_{j=i+2}^{k+1}
  (-1)^{\sumlim{\ell =i+1}{j-1} |\alpha_\ell||\alpha_i|}
  \cdot [ \ldots [\frak q^{(k)}_{n_0}(\alpha_0), \frak q_{n_1}(\alpha_1)],
  \ldots ], q_{n_{i-1}}(\alpha_{i-1})], q_{n_{i+1}}(\alpha_{i+1})], \\
 &\qquad\qquad\quad \ldots ], q_{n_{j-1}}(\alpha_{j-1})],
  [q_{n_i}'(\alpha_i), q_{n_j}(\alpha_j)]], q_{n_{j+1}}(\alpha_{j+1})],
  \ldots ], \frak q_{n_{k+1}}(\alpha_{k+1})]  \\
=& (-1)^{\sumlim{\ell =i+1}{k+1} |\alpha_\ell||\alpha_i|}
  [\Cal C(k; n_0, \alpha_0; \ldots ; n_{i-1}, \alpha_{i-1};
  n_{i+1}, \alpha_{i+1}; \ldots ; n_{k+1}, \alpha_{k+1}),
  q_{n_i}'(\alpha_i)]  \\
 &+ (-n_in_{i+1}) \cdot \Cal C(k; n_0, \alpha_0; \ldots ;
  n_{i-1}, \alpha_{i-1};
  n_i+n_{i+1}, \alpha_i\alpha_{i+1};   \\
 &\qquad\qquad\quad n_{i+2}, \alpha_{i+2};
  \ldots ; n_{k+1}, \alpha_{k+1})      \\
 &+ \sum_{j=i+2}^{k+1} (-1)^{\sumlim{\ell =i+1}{j-1}
  |\alpha_\ell||\alpha_i|} \cdot (-n_in_j) \cdot
  \Cal C(k; n_0, \alpha_0; \ldots ; n_{i-1}, \alpha_{i-1};
  n_{i+1}, \alpha_{i+1};   \\
 &\qquad\qquad\quad \ldots; n_{j-1}, \alpha_{j-1};
  n_i+n_j, \alpha_i\alpha_j; n_{j+1}, \alpha_{j+1};
  \ldots ; n_{k+1}, \alpha_{k+1})      \\
\endalign$$
where we have used Theorem 2.16 (v) in the last step.
By induction hypothesis,
$$\align
 &[ \ldots [\frak q^{(k)}_{n_0}(\alpha_0), \frak q_{n_1}(\alpha_1)],
  \ldots ], q_{n_{i-1}}(\alpha_{i-1})], q_{n_i}'(\alpha_i)],
  q_{n_{i+1}}(\alpha_{i+1})], \ldots ], \frak q_{n_{k+1}}(\alpha_{k+1})] \\
=& (-1)^{\sumlim{\ell =i+1}{k+1} |\alpha_\ell||\alpha_i|} \cdot
  [a_i \cdot \frak q_{\sum_{\ell=0}^{k+1} n_\ell - n_i}
  (\alpha_0 \ldots \alpha_{i-1}\alpha_{i+1} \ldots \alpha_{k+1}),
  q_{n_i}'(\alpha_i)]  \\
 &+ \sum_{j=i+1}^{k+1} (-n_in_j) \left \{ a_{i,j}
  \frak q_{n_0 + \ldots + n_{k+1}}(\alpha_0 \ldots \alpha_{k+1})
 +b_{i,j} \cdot \int_X(K_X \alpha_0 \ldots \alpha_{k+1})
 \cdot \Id_{\Bbb H} \right \}   \\
=&a_i \cdot \left (\sum_{\ell=0}^{k+1} n_\ell - n_i \right )n_i
  \left \{ \frak q_{n_0 + \ldots + n_{k+1}}(\alpha_0 \ldots \alpha_{k+1})
  +b_i \cdot \int_X(K_X \alpha_0 \ldots \alpha_{k+1})
  \cdot \Id_{\Bbb H} \right \}   \\
 &+ \sum_{j=i+1}^{k+1} (-n_in_j) \left \{ a_{i,j}
  \frak q_{n_0 + \ldots + n_{k+1}}(\alpha_0 \ldots \alpha_{k+1})
 +b_{i,j} \cdot \int_X(K_X \alpha_0 \ldots \alpha_{k+1})
  \cdot \Id_{\Bbb H} \right \}   \\
\endalign$$
where all the numbers $a_i, b_i, a_{i,j}, b_{i,j}$
depend only on $k, n_0, \ldots, n_k$.
Combining the above with (3.5) and (3.6), we obtain (3.2).
\endproof

\lemma{3.7} Let $k \ge 0$, $n_0, \ldots, n_k \in \Zee$,
and $\alpha_0, \ldots, \alpha_k \in H^*(X)$. Then,
$$\align
&\qquad [ \ldots [\frak q^{(k)}_{n_0}(\alpha_0), \frak q_{n_1}(\alpha_1)],
   \ldots ], \frak q_{n_k}(\alpha_k)] \\
&=\left ( (-1)^k \cdot k! \cdot n_0^k \cdot n_1 \ldots n_k \right )
   \cdot \frak q_{n_0 + \ldots + n_k}(\alpha_0 \ldots \alpha_k)
   +b \cdot \int_X(K_X \alpha_0 \ldots \alpha_k) \cdot \Id_{\Bbb H} \\
\endalign$$
where $b$ is a constant depending only on $k, n_0, \ldots, n_k$.
\endstatement
\medskip\noindent
{\it Proof.}
For simplicity, denote the constant $a$ in (3.2) by
$a(k; n_0, \ldots , n_k)$. In view of (3.2), it remains to show that
$a(k; n_0, \ldots , n_k) =
(-1)^k \cdot k! \cdot n_0^k \cdot n_1 \ldots n_k$.

This is trivially true for $k = 0$. By Theorem 2.16 (v),
this is true for $k=1$. In the following,
assuming $a(k; n_0, \ldots , n_k) =
(-1)^k \cdot k! \cdot n_0^k \cdot n_1 \ldots n_k$, we shall prove
$$a(k+1; n_0, \ldots , n_{k+1}) = (-1)^{k+1} \cdot (k+1)!
\cdot n_0^{k+1} \cdot n_1 \ldots n_{k+1}. \tag 3.8$$
Indeed, we see from the proof of Lemma 3.1 that
$a(k+1; n_0, \ldots , n_{k+1})$ is equal to
$$-\sum_{i=1}^{k+1}
\{ a(k; n_0, \ldots , n_{i-1}, n_{i+1},\ldots , n_{k+1})
\cdot \left (\sum_{\ell=0}^{k+1} n_\ell - n_i \right )n_i$$
$$+ \sum_{j=i+1}^{k+1} (-n_in_j) \cdot
a(k; n_0, \ldots , n_{i-1}, n_{i+1},\ldots ,
n_{j-1}, n_i+n_j, n_{j+1},\ldots ,n_{k+1}) \}.$$
So by using the induction hypothesis, we conclude that
$$\align
&\qquad a(k+1; n_0, \ldots , n_{k+1}) \\
&= -\sum_{i=1}^{k+1}
   \{ (-1)^k \cdot k! \cdot n_0^k \cdot n_1 \ldots n_{i-1} n_{i+1} \ldots
   n_{k+1} \cdot \left (\sum_{\ell=0}^{k+1} n_\ell - n_i \right )n_i  \\
&\quad + \sum_{j=i+1}^{k+1} (-n_in_j) \cdot
   (-1)^k \cdot k! \cdot n_0^k \cdot n_1 \ldots n_{i-1} n_{i+1} \ldots
   n_{j-1} (n_i+n_j) n_{j+1} \ldots n_{k+1} \}  \\
&= (-1)^{k+1} \cdot (k+1)!
\cdot n_0^{k+1} \cdot n_1 \ldots n_{k+1}.   \qed
\endalign$$

\lemma{3.9} Let $k \ge 0$, $n_0, \ldots, n_{k+1} \in \Zee$,
and $\alpha_0, \ldots, \alpha_{k+1} \in H^*(X)$. Then,
$$\align
&\qquad [ \ldots [\frak q^{(k)}_{n_0}(\alpha_0), \frak q_{n_1}(\alpha_1)],
   \ldots ], \frak q_{n_{k+1}}(\alpha_{k+1})] \\
&= k! \cdot n_0^k \cdot \prod_{\ell =1}^{k+1} (-n_\ell)
   \cdot \delta_{n_0 + \ldots + n_{k+1}}
   \cdot \int_X(\alpha_0 \ldots \alpha_{k+1}) \cdot \Id_{\Bbb H}. \\
\endalign$$
\endstatement
\proof
Follows from Lemma 3.7 and Theorem 2.16 (i).
\endproof

\proposition{3.10} Let $\frak f \in \End(\Bbb H)$ be an operator
of bi-degree $(0, i)$ with $\frak f' = 0$ and $\frak f^\dagger = \frak f$.
Assume that $[\frak f, \frak q_1(\alpha)]
= \sumlim{j=0}{k} c_j \frak q_1^{(j)}(\lambda_j \alpha)$
for every $\alpha \in H^*(X)$, where $k \ge 0$, $c_k = 1/k!$, and
$k, c_j \in \Bbb Q, \lambda_j \in H^{i-2j}(X)$ depend only on $\frak f$. Then,
$$\align
&\qquad [ \ldots [[\frak f, \frak q_{n_0}(\alpha_0)],
   \frak q_{n_1}(\alpha_1)], \ldots ], \frak q_{n_{k+1}}(\alpha_{k+1})] \\
&= -\prod_{\ell =0}^{k+1} (-n_\ell) \cdot \delta_{n_0 + \ldots + n_{k+1}}
   \cdot \int_X(\lambda_k \alpha_0 \ldots \alpha_{k+1}) \cdot \Id_{\Bbb H}
   \tag 3.11
\endalign$$
for all $n_0, \ldots, n_{k+1} \in \Zee$ and
all $\alpha_0, \ldots, \alpha_{k+1} \in H^*(X)$.
\endstatement
\medskip\noindent
{\it Proof.}
(i) First, we show that (3.11) holds for $n_0 \ge 0$.
We shall use induction on $n_0$. When $n_0 = 0$,
(3.11) is true since $\frak q_0(\alpha_0) = 0$ for any $\alpha_0$.
When $n_0 = 1$, we see from the assumptions on $\frak f$ that
$[\frak f, \frak q_1(\alpha_0)]
= \sum_{j=0}^k c_j \frak q_1^{(j)}(\lambda_j \alpha_0)$. So
$$\align
&\qquad [ \ldots [[\frak f, \frak q_{n_0}(\alpha_0)],
   \frak q_{n_1}(\alpha_1)], \ldots ], \frak q_{n_{k+1}}(\alpha_{k+1})] \\
&= \sum_{j=0}^k c_j
   [ \ldots [ \frak q_1^{(j)}(\lambda_j \alpha_0), \frak q_{n_1}(\alpha_1)],
   \ldots ], \frak q_{n_{k+1}}(\alpha_{k+1})].    \tag 3.12
\endalign$$
By Lemma 3.9,
$[ \ldots [ \frak q_1^{(j)}(\lambda_j \alpha_0), \frak q_{n_1}(\alpha_1)],
\ldots ], \frak q_{n_{k+1}}(\alpha_{k+1})]=0$ if $0 \le j < k$. Also,
$$\align
&\qquad [ \ldots [ \frak q_1^{(k)}(\lambda_k \alpha_0),
   \frak q_{n_1}(\alpha_1)], \ldots ], \frak q_{n_{k+1}}(\alpha_{k+1})]  \\
&= -k! \cdot \prod_{\ell =0}^{k+1} (-n_\ell) \cdot
   \delta_{n_0 + \ldots + n_{k+1}} \cdot
   \int_X(\lambda_k \alpha_0 \ldots \alpha_{k+1})
   \cdot \Id_{\Bbb H}  \tag 3.13
\endalign$$
by Lemma 3.9 again. It follows from (3.12) and (3.13) that
(3.11) holds for $n_0=1$.

Next, assuming that (3.11) is true for some positive integer $n_0$,
we shall prove that (3.11) still holds if $n_0$ is replaced by $(n_0+1)$.
Note from Theorem 2.16 (v) that $\frak q_{n_0+1}(\alpha_0) = - {1/n_0}
\cdot [\frak q_1'(1_X), \frak q_{n_0}(\alpha_0)]$. Thus, we obtain
$$\align
 [\frak f, \frak q_{n_0+1}(\alpha_0)]
&= - {1 \over n_0} \cdot
 [\frak f, [\frak q_1'(1_X), \frak q_{n_0}(\alpha_0)]] \\
&= - {1 \over n_0} \cdot
 \{ [[\frak f, \frak q_1'(1_X)], \frak q_{n_0}(\alpha_0)]
 + [\frak q_1'(1_X), [\frak f, \frak q_{n_0}(\alpha_0)]]\}.  \\
\endalign$$
Since $\frak f' =0$,
$[\frak f, \frak q_1'(1_X)] = [\frak f, \frak q_1(1_X)]'
= \sum_{j=0}^k c_j \frak q_1^{(j+1)}(\lambda_j)$. Therefore, we have
$$[\frak f, \frak q_{n_0+1}(\alpha_0)]
= - {1 \over n_0} \left \{ \sum_{j=0}^k c_j
    [\frak q_1^{(j+1)}(\lambda_j), \frak q_{n_0}(\alpha_0)]
- [[\frak f, \frak q_{n_0}(\alpha_0)], \frak q_1'(1_X)] \right \}.
\tag 3.14$$
Note that $[ \ldots [ \frak q_1^{(j+1)}(\lambda_j), \frak q_{n_0}(\alpha_0)],
\ldots ], \frak q_{n_{k+1}}(\alpha_{k+1})]$ is equal to $0$ if $0 \le j < k$,
and $(k+1)! \cdot \prod_{\ell =0}^{k+1} (-n_\ell) \cdot
\delta_{1+n_0 + \ldots + n_{k+1}} \cdot
\int_X(\lambda_k \alpha_0 \ldots \alpha_{k+1})
\cdot \Id_{\Bbb H}$ if $j = k$. So
$$\align
&\qquad [ \ldots [[\frak f, \frak q_{n_0+1}(\alpha_0)],
   \frak q_{n_1}(\alpha_1)], \ldots ], \frak q_{n_{k+1}}(\alpha_{k+1})] \\
&= -{1 \over n_0} \{ (-1)^k(k+1) \cdot \prod_{\ell =0}^{k+1} n_\ell \cdot
   \delta_{1+n_0 + \ldots + n_{k+1}} \cdot
   \int_X(\lambda_k \alpha_0 \ldots \alpha_{k+1}) \cdot \Id_{\Bbb H}   \\
&\qquad - [ \ldots [[\frak f, \frak q_{n_0}(\alpha_0)], \frak q_1'(1_X)],
   \frak q_{n_1}(\alpha_1)], \ldots ], \frak q_{n_{k+1}}(\alpha_{k+1})] \} \\
&= -{1 \over n_0} \{ (-1)^k(k+1) \cdot \prod_{\ell =0}^{k+1} n_\ell \cdot
   \delta_{1+n_0 + \ldots + n_{k+1}} \cdot
   \int_X(\lambda_k \alpha_0 \ldots \alpha_{k+1}) \cdot \Id_{\Bbb H}   \\
&\qquad - [ \ldots [[\frak f, \frak q_{n_0}(\alpha_0)],
   \frak q_{n_1}(\alpha_1)],
   \ldots ], \frak q_{n_{k+1}}(\alpha_{k+1})], \frak q_1'(1_X)]  \\
&\qquad - \sum_{j=1}^{k+1} [ \ldots [[\frak f, \frak q_{n_0}(\alpha_0)],
   \ldots ], \frak q_{n_{j-1}}(\alpha_{j-1})],
   [q_1'(1_X),\frak q_{n_j}(\alpha_j)]], \frak q_{n_{j+1}}(\alpha_{j+1})],
   \ldots    \\
&\qquad\qquad\qquad \ldots], \frak q_{n_{k+1}}(\alpha_{k+1})] \}. \\
\endalign$$
By induction hypothesis,
$[ \ldots [[\frak f, \frak q_{n_0}(\alpha_0)], \frak q_{n_1}(\alpha_1)],
\ldots ], \frak q_{n_{k+1}}(\alpha_{k+1})], \frak q_1'(1_X)] = 0$.
By Theorem 2.16 (v) and induction hypothesis, we get
$$\align
&\qquad [ \ldots [[\frak f, \frak q_{n_0}(\alpha_0)],
   \ldots ], \frak q_{n_{j-1}}(\alpha_{j-1})],
   [q_1'(1_X),\frak q_{n_j}(\alpha_j)]], \frak q_{n_{j+1}}(\alpha_{j+1})],
   \ldots    \\
&\qquad\qquad\qquad \ldots], \frak q_{n_{k+1}}(\alpha_{k+1})] \\
&= (-n_j) \cdot [ \ldots [[\frak f, \frak q_{n_0}(\alpha_0)],
   \ldots ], \frak q_{n_{j-1}}(\alpha_{j-1})],
   \frak q_{n_j+1}(\alpha_j)], \frak q_{n_{j+1}}(\alpha_{j+1})],
   \ldots    \\
&\qquad\qquad\qquad \ldots], \frak q_{n_{k+1}}(\alpha_{k+1})]   \\
&= (1+n_j) \cdot \prod_{\ell =0}^{k+1} (-n_\ell) \cdot
   \delta_{1+n_0 + \ldots + n_{k+1}} \cdot
   \int_X(\lambda_k \alpha_0 \ldots \alpha_{k+1}) \cdot \Id_{\Bbb H}. \\
\endalign$$
It follows that
$[ \ldots [[\frak f, \frak q_{n_0+1}(\alpha_0)], \frak q_{n_1}(\alpha_1)],
\ldots ], \frak q_{n_{k+1}}(\alpha_{k+1})]$ is equal to
$$-(-(1+n_0)) \prod_{\ell =1}^{k+1} (-n_\ell) \cdot
\delta_{1+n_0 + \ldots + n_{k+1}} \cdot
\int_X(\lambda_k \alpha_0 \ldots \alpha_{k+1}) \cdot \Id_{\Bbb H},
\tag 3.15$$
i.e., (3.11) is still true if $n_0$ is replaced by $(n_0+1)$.

(ii) Now we show that (3.11) holds for $n_0 < 0$. Note that
$\frak q_{-n}(\alpha) = (-1)^n \frak q_{n}(\alpha)^\dagger$,
$(\frak g_1 \frak g_2)^\dagger = (-1)^{m_1m_2} \cdot
(\frak g_2)^\dagger (\frak g_1)^\dagger$
and $[\frak g_1, \frak g_2]^\dagger
= -[(\frak g_1)^\dagger, (\frak g_2)^\dagger]$
for $\frak g_1, \frak g_2 \in \End(\Bbb H)$ of bi-degrees
$(\ell_1, m_1)$, $(\ell_2, m_2)$ respectively.
Since $\frak f^\dagger = \frak f$ by assumption,
$$\align
&\qquad [ \ldots [[\frak f, \frak q_{n_0}(\alpha_0)],
   \frak q_{n_1}(\alpha_1)], \ldots ], \frak q_{n_{k+1}}(\alpha_{k+1})] \\
&= (-1)^{(1+n_0)+\ldots+(1+n_{k+1})}
   \cdot [ \ldots [[\frak f, \frak q_{-n_0}(\alpha_0)],
   \frak q_{-n_1}(\alpha_1)],
   \ldots ], \frak q_{-n_{k+1}}(\alpha_{k+1})]^\dagger.
\endalign$$
By what we have proved in (i) for the positive integer $(-n_0)$,
we have
$$\align
&\qquad [ \ldots [[\frak f, \frak q_{n_0}(\alpha_0)],
   \frak q_{n_1}(\alpha_1)], \ldots ], \frak q_{n_{k+1}}(\alpha_{k+1})] \\
&=(-1)^{n_0+\ldots+n_{k+1}+k} \cdot
   \left ( -\prod_{\ell =0}^{k+1} n_\ell \cdot
   \delta_{-n_0 - \ldots - n_{k+1}} \cdot
   \int_X(\lambda_k \alpha_0 \ldots \alpha_{k+1})
   \cdot \Id_{\Bbb H} \right )^\dagger \\
&=-\prod_{\ell =0}^{k+1} (-n_\ell) \cdot
   \delta_{n_0 + \ldots + n_{k+1}} \cdot
   \int_X(\lambda_k \alpha_0 \ldots \alpha_{k+1})
   \cdot \Id_{\Bbb H}. \qed
\endalign$$

Next, we shall define the transfer property for certain commutators,
and verify that the commutator
$[\frak q^{(k)}_{n_0}(\alpha_0), \frak q_{n_1}(\alpha_1)]$
satisfies the transfer property.

\medskip\noindent
{\bf Definition 3.16.} Let $\frak A(\alpha), \frak B(\beta)
\in \End(\Bbb H)$ be two series of operators depending linearly on
$\alpha, \beta \in H^*(X)$. Then, the commutator
$[\frak A(\alpha), \frak B(\beta)]$ satisfies {\it the transfer property}
if $[\frak A(\alpha), \frak B(\beta)]
= [\frak A(1_X), \frak B(\alpha \beta)]
= [\frak A(\alpha \beta), \frak B(1_X)]$
for all $\alpha, \beta$.

\proposition{3.17} Let $k$ be a nonnegative integer.
Let $n_0, n_1 \in \Zee$, and $\alpha_0, \alpha_1 \in H^*(X)$.
Then, the commutator $[\frak q^{(k)}_{n_0}(\alpha_0),
\frak q_{n_1}(\alpha_1)]$ satisfies the transfer property:
$$[\frak q^{(k)}_{n_0}(\alpha_0), \frak q_{n_1}(\alpha_1)]
= [\frak q^{(k)}_{n_0}(1_X), \frak q_{n_1}(\alpha_0 \alpha_1)]
= [\frak q^{(k)}_{n_0}(\alpha_0 \alpha_1), \frak q_{n_1}(1_X)]. \tag 3.18$$
\endstatement
\proof
By Theorem 2.16 (i) and (v), (3.18) is true for $k = 0, 1$.
In the following, we assume that $k \ge 2$. In addition, we may assume
that $(\alpha_0, \alpha_1) \ne (1_X, 1_X)$. Let $\frak e$ be the difference
of any of the two commutators in (3.18). By Lemma 3.1, we have
$$[ \ldots [\frak e, \frak q_{n_2}(\alpha_2)],
\ldots ], \frak q_{n_{k-1}}(\alpha_{k-1})], \frak q_{n_k}(\alpha_k)] = 0
\tag 3.19$$
for all $n_2, \ldots, n_k \in \Zee$ and
all $\alpha_2, \ldots, \alpha_k \in H^*(X)$.
Since $\Bbb H$ is irreducible (see the last paragraph in section 2),
we see from Schur's lemma that
$$\frak e_{k-1} {\overset \text{def} \to =}
[ \ldots [\frak e, \frak q_{n_2}(\alpha_2)],
\ldots ], \frak q_{n_{k-1}}(\alpha_{k-1})]$$
must be a scalar multiple of the identity operator. Now,
the bi-degree of $\frak e_{k-1}$ is
$$\left (\sum_{j=0}^{k-1} n_j, 2k+ \sum_{j=0}^{k-1} (2n_j-2+|\alpha_j|)
\right ).$$
So when $\sum_{j=0}^{k-1} n_j \ne 0$, the bi-degree of $\frak e_{k-1}$ is
nontrivial. When $\sum_{j=0}^{k-1} n_j = 0$,
the bi-degree of $\frak e_{k-1}$ is
$(0, \sum_{j=0}^{k-1} |\alpha_j|)$ which again is nontrivial
since we have assumed $(\alpha_0, \alpha_1) \ne (1_X, 1_X)$.
This forces $\frak e_{k-1} = 0$, i.e., we have
$$[ \ldots [\frak e, \frak q_{n_2}(\alpha_2)], \ldots ],
\frak q_{n_{k-2}}(\alpha_{k-2})], \frak q_{n_{k-1}}(\alpha_{k-1})] = 0$$
for all the integers $n_2, \ldots, n_{k-1} \in \Zee$ and
all the cohomology classes $\alpha_2, \ldots, \alpha_{k-1} \in H^*(X)$.
Repeating the above process, we conclude that $\frak e = 0$.
\endproof

\section{4. The operators $W_n^k$ as the leading terms}

In this section, we shall define certain operators $W_n^k(\alpha)$
and determine their commutation relations with the Heisenberg generator
$\frak q_m(\beta)$. Our main result (see Theorem 4.12 below)
says that if $\frak f \in \End(\Bbb H)$
satisfies the same conditions as in Proposition 3.10,
then the leading term of $\frak f$ is equal to $-W_0^{k+2}(\lambda_k)$.

First, we recall the normally ordered product
$: \cdot :$ from \cite{Bor, FLM, Kac}. Let
$$a(z) = \sum_{n \in \Zee} a_{(n)} z^{n-\Delta}$$
be a vertex operator of conformal weight $\Delta$,
that is, a generating function in a formal variable $z$
with $a_{(n)} \in \End(\Bbb H)$ of bi-degree $(n, *)$. Put
$$a_+(z)= \sum_{n >0} a_{(n)} z^{n-\Delta} \qquad \text{and} \qquad
a_-(z) = \sum_{n \le 0} a_{(n)} z^{n-\Delta}$$
(note that our sign convention on vertex operators throughout this
paper differs from the standard one used in the vertex algebra
literature \cite{Bor, FLM, Kac}).
If $b(z)$ is another vertex operator, we define a new vertex operator,
which is called {\it the normally ordered product} of $a(z)$ and $b(z)$,
to be:
$$ :a(z)b(z): = a_+(z) b(z) + (-1)^{ab} b(z) a_-(z)  \tag 4.1$$
where $(-1)^{ab}$ is $-1$ if both $a(z)$ and $b(z)$ are odd
and $1$ otherwise. Inductively we can define
the normally ordered product of $k$ vertex operators from right to left by
$$:a_1(z) a_2(z) \cdots a_k(z): \, = \,
: a_1(z) (:a_2(z) \cdots a_k(z):):.  \tag 4.2$$

\medskip\noindent
{\bf Definition 4.3.} (i) For $\alpha \in H^*(X)$,
we define a vertex operator $\alpha(z)$ by putting
$$ \alpha(z) = \sum_{n \in \Zee} \frak q_n(\alpha) z^{n-1}; \tag 4.4$$

(ii) Let $k \ge 1$, $n \in \Zee$, and $\alpha \in H^*(X)$.
Let $\tau_{k*}: H^*(X) \to H^*(X^k)$ be the linear map induced by
the diagonal embedding $\tau_k: X \to X^k$, and
$$\tau_{k*}\alpha = \sum_j \alpha_{j,1} \otimes \ldots \otimes
\alpha_{j, k}$$
via the K\"unneth decomposition of $H^*(X^k)$.
We define the operator $W^k_n(\alpha) \in \End(\Bbb H)$
to be the coefficient of $z^{n-k}$ in the vertex operator
$$ \frac1{k!} \cdot (\tau_{k*}\alpha) (z) \overset \text{def} \to =
\frac1{k!} \cdot \sum_j : \alpha_{j, 1} (z) \cdots \alpha_{j, k}(z):.
\tag 4.5$$

Note that $W_n^k(\alpha)$ is a homogeneous linear operator of bi-degree
$$(n, 2n+2k-4+|\alpha|).$$
Also, $W_n^1(\alpha)$ coincides with the Heisenberg generator
$\frak q_n(\alpha)$, and $W_n^2(\alpha)$ coincides with
the Virasoro generator $\frak L_n(\alpha)$.
The next lemma generalizes Theorem 2.16 (ii) and indicates that
the commutator $[W_n^k(\alpha), \frak q_m(\beta)]$
satisfies the transfer property.

\lemma{4.6} Let $k \ge 2$. Let $n, m \in \Zee$ and
$\alpha, \beta \in H^*(X)$. Then,
$$[W_n^k(\alpha), \frak q_m(\beta)]
= (-m) \cdot W_{n+m}^{k-1}(\alpha \beta).  \tag 4.7$$
\endstatement
\medskip\noindent
{\it Proof.}
First of all, we rewrite the commutation relation Theorem 2.16 (i) as
$$[{\alpha} (z), \beta(w)] = \int_X (\alpha \beta) \cdot
\sum_{n \in \Zee} (n \cdot z^{n-1}w^{-n-1}).  $$
Assume that $\tau_{k*}\alpha = \sum_j \alpha_{j,1} \otimes \ldots
\otimes \alpha_{j, k} \in H^*(X^k)$. Then, we have
$$\align
 &\qquad [(\tau_{k*}\alpha)(z), \frak q_m(\beta) ]   \\
 & = \Res{w} w^{-m} [(\tau_{k*}\alpha)(z), \beta(w) ] \\
 & =  \Res{w} w^{-m} \sum_j
    [: \alpha_{j,1}(z) \cdots \alpha_{j,k}(z): , \beta(w) ]  \\
 & =  \Res{w} w^{-m} \sum_{s=1}^k \sum_j
    : \alpha_{j,1}(z) \cdots \widehat{\alpha}_{j,s}(z)
      \cdots \alpha_{j,k}(z): \cdot    \\
 &\qquad \cdot (-1)^{|\beta| \sum_{\ell=s+1}^k |\alpha_{j,\ell}|}
      \int_X \alpha_{j,s} \beta \sum_{n \in \Zee} n z^{n-1}w^{-n-1}\\
 & =  (-m z^{-m-1}) \cdot\sum_{s=1}^k \sum_j
    : \alpha_{j,1}(z) \cdots \widehat{\alpha}_{j,s}(z)
      \cdots \alpha_{j,k}(z): \cdot    \\
 &\qquad \cdot (-1)^{|\beta| \sum_{\ell=s+1}^k |\alpha_{j,\ell}|}
      \int_X \alpha_{j,s}\beta. \tag 4.8  \\
\endalign$$
To simplify (4.8), we fix $s$ satisfying $1 \le s \le k$.
Let $p_s$ (respectively, $\pi_s$) be
the projection of $X^k$ to the $s$-th factor (respectively,
to the product of the remaining $(k-1)$ factors).
So we have a commutative diagram of morphisms:
$$\matrix
X &\overset \tau_k \to \longrightarrow &X^k \quad
           &\overset p_s \to \longrightarrow &X \\
  &\searrow^{\tau_{k-1}} &\downarrow{\pi_s}& &  \\
  &                         &X^{k-1} \quad & &  \\
\endmatrix  $$
Now, by the projection formula, we conclude that
$$\align
&\quad (\tau_{k-1})_* (\alpha \beta)
 = (\pi_s \circ \tau_k)_* (\alpha \beta) \\
&= \pi_{s*}(\tau_{k*}(\alpha \cdot \tau_k^* p_s^* \beta ))
 = \pi_{s*}( \tau_{k*} \alpha \cdot p_s^* \beta ) \\
&= \sum_j \alpha_{j,1} \otimes \cdots \otimes \widehat{\alpha}_{j,s}
 \otimes \cdots \otimes \alpha_{j, k} \cdot
 (-1)^{|\beta| \sum_{\ell=s+1}^k |\alpha_{j,\ell}|}
 \int_X \alpha_{j,s}\beta.
\endalign$$
Combining this with (4.8), we conclude that
$$[(\tau_{k*}\alpha) (z), \frak q_m(\beta) ] = (-m z^{-m-1})
\cdot k \cdot ((\tau_{k-1})_* (\alpha \beta))(z). \tag 4.9$$
Now comparing the coefficients of $z^{n-k}$ on both sides of (4.9),
we obtain (4.7). \qed

\proposition{4.10} Let $k \ge 1, n_0, \ldots, n_k \in \Zee$,
and $\alpha_0, \ldots, \alpha_k \in H^*(X)$. Then,
\par
{\rm (i)} $[ \ldots [W^k_{n_0}(\alpha_0), \frak q_{n_1}(\alpha_1)],
   \ldots ], \frak q_{n_{k-1}}(\alpha_{k-1})]$ is equal to
$$\prod_{\ell =1}^{k-1} (-n_\ell) \cdot
   \frak q_{n_0 + \ldots + n_{k-1}}(\alpha_0 \ldots \alpha_{k-1});$$
\par
{\rm (ii)} $[ \ldots [W^k_{n_0}(\alpha_0), \frak q_{n_1}(\alpha_1)],
   \ldots ], \frak q_{n_k}(\alpha_k)]$ is equal to
$$\prod_{\ell =1}^{k} (-n_\ell)
   \cdot \delta_{n_0 + \ldots + n_{k}}
   \cdot \int_X(\alpha_0 \ldots \alpha_{k}) \cdot \Id_{\Bbb H}.$$
\endstatement
\proof
Applying Lemma 4.6 repeatedly, we see that
$$\align
&\qquad [ \ldots [W^k_{n_0}(\alpha_0), \frak q_{n_1}(\alpha_1)],
   \ldots ], \frak q_{n_{k-1}}(\alpha_{k-1})]         \\
&= \prod_{\ell =1}^{k-1} (-n_\ell) \cdot
   W^1_{n_0 + \ldots + n_{k-1}}(\alpha_0 \ldots \alpha_{k-1})  \\
&= \prod_{\ell =1}^{k-1} (-n_\ell) \cdot
   \frak q_{n_0 + \ldots + n_{k-1}}(\alpha_0 \ldots \alpha_{k-1}).  \\
\endalign$$
This proves (i). Now (ii) follows from (i) and Theorem 2.16 (i).
\endproof

Let $\frak f \in \End(\Bbb H)$ be a linear operator.
We write the operator $\frak f$
as a (possibly infinite) linear combination of
finite products of Heisenberg generators:
$$\frak q_{m_1}(\beta_1) \cdots \frak q_{m_i}(\beta_i)  \tag 4.11$$
where those $\frak q_{m_j}(\beta_j)$ with $m_j <0$ are put to the right.
In other words, we regard $\frak f$ as an element in the completion
of the universal enveloping algebra of the Heisenberg algebra.
We can always do so because $\Bbb H$ is irreducible as
a representation of the Heisenberg algebra.
Assume that the lengths $i$ of all the product terms (4.11) appearing in
$\frak f$ have a common upper bound (this is the case for all
the operators $\frak f$ considered below).
Then we define the {\it leading term} of $\frak f$ to be the sum of
those products (4.11) in $\frak f$ such that $i$ is the largest.
For example, each product in the operator $W^{k}_n(\alpha)$
has $k$ factors of $\frak q$'s by Definition 4.3 (ii).
So the leading term of $W^{k}_n(\alpha)$ is itself.
Our next result says that if $\frak f \in \End(\Bbb H)$
satisfies the same assumptions as in Proposition 3.10,
then the leading term of $\frak f$ is $-W_0^{k+2}(\lambda_k)$.

\theorem{4.12} Let $\frak f \in \End(\Bbb H)$ be of bi-degree $(0, i)$
with $\frak f' = 0$ and $\frak f^\dagger = \frak f$.
Assume $[\frak f, \frak q_1(\alpha)]
= \sumlim{j=0}{k} c_j \frak q_1^{(j)}(\lambda_j \alpha)$
for every $\alpha \in H^*(X)$, where $k \ge 0$, $c_k = 1/k!$,
and $k, c_j \in \Bbb Q, \lambda_j \in H^{i-2j}(X)$ depend only on $\frak f$.
Put $\epsilon(\frak f) = \frak f + W_0^{k+2}(\lambda_k)$. Then,
\par
{\rm (i)} for all $n_1, \ldots, n_{k+2} \in \Zee$ and
all $\alpha_1, \ldots, \alpha_{k+2} \in H^*(X)$,
$$[ \ldots [\epsilon(\frak f), \frak q_{n_1}(\alpha_1)],
\ldots ], \frak q_{n_{k+2}}(\alpha_{k+2})] = 0;$$
\par
{\rm (ii)} for all $n_1, \ldots, n_{k+1} \in \Zee$ with
$\sum_{j = 1}^{k+1} n_j \ne 0$ and
all $\alpha_1, \ldots, \alpha_{k+1} \in H^*(X)$,
$$[ \ldots [\epsilon(\frak f), \frak q_{n_1}(\alpha_1)],
\ldots ], \frak q_{n_{k+1}}(\alpha_{k+1})] = 0;$$
\par
{\rm (iii)} the leading term of $\frak f$ is $-W_0^{k+2}(\lambda_k)$.
\endstatement
\proof
(i) Follows from Proposition 3.10 and Proposition 4.10 (ii).

(ii) Denote $[ \ldots [\epsilon(\frak f), \frak q_{n_1}(\alpha_1)],
\ldots ], \frak q_{n_{k+1}}(\alpha_{k+1})]$ by $\frak g$. Then,
we see from (i) that $[\frak g, \frak q_{n_{k+2}}(\alpha_{k+2})] = 0$
for all $n_{k+2} \in \Zee$ and all $\alpha_{k+2} \in H^*(X)$.
Since $\Bbb H$ is irreducible, $\frak g$ must be a scalar.
Now, the bi-degree of $\frak g$ is equal to
$$\left (\sum_{j=1}^{k+1} n_j, i + \sum_{j=1}^{k+1} (2n_j-2+|\alpha_j|)
\right )$$
which is nontrivial since $\sum_{j=1}^{k+1} n_j \ne 0$.
Therefore, $\frak g =0$.

(iii) As in (4.11), we write the operator $\frak f$ as
a (possibly infinite)
linear combination of finite products of the Heisenberg generators:
$$\frak q_{m_1}(\beta_1) \cdots \frak q_{m_i}(\beta_i)  \tag 4.13$$
where those $\frak q_{m_j}(\beta_j)$ with $m_j <0$ are put to the right.
Note that
$$\frak f=\epsilon(\frak f)-W_0^{k+2}(\lambda_k).$$
Now (i) says that $\epsilon(\frak f)$ is a (possibly infinite)
linear combination of finite products of
at most $(k+1)$ Heisenberg generators.
Also, from the discussion in the paragraph preceding Theorem 4.12,
we see that each product in the operator $W_0^{k+2}(\lambda_k)$
has $(k+2)$ factors of $\frak q$'s. Therefore,
the leading term of $\frak f$ is $-W_0^{k+2}(\lambda_k)$.
\endproof

\section{5. The cohomology ring structure}

In this section, we study   the cohomology ring structure
of
$\Bbb H_n = H^*(\hil{n})$. In particular, we find the ring
generators  of  $H^*(\hil{n})$. The basic idea is to introduce certain
operator
$\frak G(\gamma) \in \End(\Bbb H)$ for $\gamma \in H^*(X)$,
which generalizes the considerations in \cite{Leh}.
We show that $\frak G(\gamma)$ satisfies the assumptions in Theorem 4.12.
So the leading term of $\frak G(\gamma)$ is related to
the operator $W_n^k(\alpha)$ from Definition 4.3 (ii).
Then an inductive procedure proves our main result.

\medskip\noindent
{\bf Definition 5.1.} (i) For $\gamma \in H^*(X)$ and $n \ge 0$, define
$$G(\gamma, n) =
p_{1*}(\ch(\Cal O_{\Cal Z_n}) \cdot p_2^*\td(X)
\cdot p_2^*\gamma) \in H^*(\hil{n})  \tag 5.2$$
where $\ch(\Cal O_{\Cal Z_n})$ is the Chern character of the sheaf
$\Cal O_{\Cal Z_n}$, $\td(X)$ is the Todd class of $X$,
and $p_1, p_2$ are the natural projections
of $\hil{n} \times X$ to $\hil{n}, X$ respectively.
We define the linear operator $\frak G(\gamma) \in \End(\Bbb H)$ by putting
$$\frak G(\gamma) = \bigopluslim{n \ge 0}{} G(\gamma, n)$$
where $G(\gamma, n)$ acts on $\Bbb H_n = H^*(\hil{n})$ by
the cup product;

(ii) For $i \in \Zee$ and $\gamma \in H^s(X)$,
define $G_i(\gamma, n)$ to be the component
of $G(\gamma, n)$ in $H^{s+2i}(\hil{n})$. We define the operator
$\frak G_i(\gamma) \in \End(\Bbb H)$ by
$$\frak G_i(\gamma) = \bigopluslim{n \ge 0}{} G_i(\gamma, n)$$
where again $G_i(\gamma, n)$ acts on $\Bbb H_n = H^*(\hil{n})$ by
the cup product.

Notice that $\frak G_i(\gamma) \in \End(\Bbb H)$ is
homogeneous of bi-degree $(0, |\gamma|+2i)$. Moreover,
$$\frak G_i(\gamma)' = 0 \qquad \text{and} \qquad
\frak G_i(\gamma)^\dagger = \frak G_i(\gamma).  \tag 5.3$$
By comparing the degrees on both sides of (5.2), we see that
$$G(\gamma, n) = \sum_{i \in \Zee} G_i(\gamma, n).$$
Therefore, by the definition of the operators $\frak G(\gamma)$ and
$\frak G_i(\gamma)$, we have
$$\frak G(\gamma) = \sum_{i \in \Zee} \frak G_i(\gamma). \tag 5.4$$

When $\gamma = 1_X$, we have $G(1_X, n) =
p_{1*}(\ch(\Cal O_{\Cal Z_n}) \cdot p_2^*\td(X))$.
So we see from the Grothendieck-Riemann-Roch Theorem \cite{Har} that
$$G(1_X, n) = \ch(p_{1*}\Cal O_{\Cal Z_n})  \tag 5.5$$
where we have made no distinction between an algebraic cycle and
its corresponding cohomology class. In particular,
we have the following two formulas:
$$\align
 \frak G_0(1_X)
&= \bigopluslim{n}{} G_0(1_X, n)
   = \bigopluslim{n}{} (n \cdot \Id_{\hil{n}})
   = \frak L_0(-1_X),   \tag 5.6   \\
 \frak G_1(1_X)
&= \bigopluslim{n}{} G_1(1_X, n)
   = \bigopluslim{n}{} c_1(p_{1*}\Cal O_{\Cal Z_n})
   = \frak d.           \tag 5.7
\endalign$$

Our next lemma and its proof are parallel to the Theorem 4.2
and its proof in \cite{Leh}. Together with (5.3), this lemma
enables us to apply Theorem 4.12.

\lemma{5.8} Let $\gamma, \alpha \in H^*(X)$. Then, we have
$$[\frak G(\gamma), \frak q_1(\alpha)]
= \text{\rm exp}(\ad(\frak d))(\frak q_1(\gamma \alpha));$$
or equivalently, using the component $\frak G_k(\gamma)$ of
$\frak G(\gamma)$ (see (5.4)), we obtain
$$[\frak G_k(\gamma), \frak q_1(\alpha)] =
{1 \over k!} \cdot \frak q_1^{(k)}(\gamma \alpha). \tag 5.9$$
\endstatement
\proof
Recall the standard diagram (2.2) and the exact sequence (2.3). We have
$$\psi_X^*\ch(\Cal O_{\Cal Z_n}) -
\phi_X^*\ch(\Cal O_{\Cal Z_{n-1}}) = \rho_X^*\ch(\Cal O_{\Delta_X})
\cdot p_1^*\text{exp}(-[E_n])  \tag 5.10$$
where we have used the fact that
$\ch(\rho_X^*\Cal O_{\Delta_X} \otimes
  p_1^*\Cal O_{\hil{n,n-1}}(-E_n))
= \rho_X^*\ch(\Cal O_{\Delta_X}) \cdot
  p_1^*\ch(\Cal O_{\hil{n,n-1}}(-E_n))
= \rho_X^*\ch(\Cal O_{\Delta_X}) \cdot p_1^*\text{exp}(-[E_n])$.

\medskip\noindent
{\bf Claim.} $\psi^*G(\gamma, n) = \phi^*G(\gamma, n-1)
+ \rho^*\gamma \cdot \text{exp}(-[E_n])$.
\medskip\noindent
{\it Proof.}
For $i \ge 1$, let $p_{i,1}$ and $p_{i,2}$ be the projections of
$\hil{i} \times X$ to $\hil{i}$ and $X$ respectively.
>From the definition of $G(\gamma, n)$, we see that
$$\align
 \psi^*G(\gamma, n)
&= \psi^*(p_{n,1})_*(\ch(\Cal O_{\Cal Z_n})
  \cdot (p_{n,2})^*\td(X) \cdot (p_{n,2})^*\gamma)   \\
&= p_{1*}\psi_X^*(\ch(\Cal O_{\Cal Z_n})
  \cdot (p_{n,2})^*\td(X) \cdot (p_{n,2})^*\gamma)   \\
&= p_{1*}(\psi_X^*\ch(\Cal O_{\Cal Z_n})
  \cdot p_2^*\td(X) \cdot p_2^*\gamma)   \\
\endalign$$
where $p_{2}$ is the projection of $\hil{n, n-1} \times X$ to $X$.
Similarly,
$$\phi^*G(\gamma, n-1) = p_{1*}
(\phi_X^*\ch(\Cal O_{\Cal Z_{n-1}}) \cdot p_2^*\td(X)
\cdot p_2^*\gamma).$$
So combining these with the formula (5.10) above, we conclude that
$$\align
&\qquad \psi^*G(\gamma, n) - \phi^*G(\gamma, n-1) \\
&= p_{1*}((\psi_X^*\ch(\Cal O_{\Cal Z_n}) -
\phi_X^*\ch(\Cal O_{\Cal Z_{n-1}})) \cdot p_2^*\td(X)
\cdot p_2^*\gamma) \\
&= p_{1*}(\rho_X^*\ch(\Cal O_{\Delta_X})
\cdot p_1^*\text{exp}(-[E_n]) \cdot p_2^*\td(X)
\cdot p_2^*\gamma) \\
&= p_{1*}(\rho_X^*\ch(\Cal O_{\Delta_X})
\cdot p_2^*\td(X) \cdot p_2^*\gamma) \cdot \text{exp}(-[E_n]) \\
&= p_{1*}\rho_X^*(\ch(\Cal O_{\Delta_X})
\cdot (p_{1,2})^*\td(X) \cdot (p_{1,2})^*\gamma) \cdot \text{exp}(-[E_n]) \\
&= \rho^* (p_{1,1})_*(\ch(\Cal O_{\Delta_X})
\cdot (p_{1,2})^*\td(X) \cdot (p_{1,2})^*\gamma) \cdot \text{exp}(-[E_n]). \\
\endalign$$
Thus, it remains to show that over $X \times X$, we have
$$(p_{1,1})_*(\ch(\Cal O_{\Delta_X})
\cdot (p_{1,2})^*\td(X) \cdot (p_{1,2})^*\gamma) = \gamma.$$
Indeed, applying the Grothendieck-Riemann-Roch Theorem to
the diagonal embedding $\tau_2: X \to X \times X$,
we get $\ch(\Cal O_{\Delta_X}) \cdot (p_{1,2})^*\td(X)
= \tau_{2*}[X]$. So
$$\align
&\qquad (p_{1,1})_*(\ch(\Cal O_{\Delta_X})
  \cdot (p_{1,2})^*\td(X) \cdot (p_{1,2})^*\gamma) \\
&= (p_{1,1})_*(\tau_{2*}[X] \cdot (p_{1,2})^*\gamma)
 = \gamma. \qed    \\
\endalign$$

We continue the proof of the lemma.
By (2.11), for any $a \in H^*(\hil{n-1})$,
$$\frak q_1(\alpha)(a) = \w p_{1*}([Q^{[n, n-1]}]
\cdot \w \rho^*\alpha \cdot \w p_2^*a).$$
Define $\iota: \hil{n, n-1} \to \hil{n} \times X \times \hil{n-1}$
by putting $\iota(\xi, \eta) = (\xi, \rho(\xi, \eta), \eta)$.
Then, $\iota$ is an embedding. Moreover, $\iota_*[X^{[n, n-1]}]
= [Q^{[n, n-1]}]$. Thus, we have
$$\align
 \frak q_1(\alpha)(a)
&= \w p_{1*}(\iota_*[X^{[n, n-1]}] \cdot
   \w \rho^*\alpha \cdot \w p_2^*a)          \\
&= \w p_{1*} \iota_*([X^{[n, n-1]}] \cdot (\iota^* \circ \w \rho^*)\alpha
   \cdot (\iota^* \circ \w p_2^*)a)          \\
&= \psi_*([X^{[n, n-1]}] \cdot \rho^*\alpha \cdot \phi^*a).
\endalign$$
Combining this with the above Claim, we conclude that
$$\align
 \frak G(\gamma) \frak q_1(\alpha)(a)
&= G(\gamma, n) \cdot \psi_*([X^{[n, n-1]}] \cdot
   \rho^*\alpha \cdot \phi^*a) \\
&= \psi_*(\psi^*G(\gamma, n) \cdot [X^{[n, n-1]}] \cdot
   \rho^*\alpha \cdot \phi^*a) \\
&= \psi_*(\phi^*G(\gamma, n-1) \cdot [X^{[n, n-1]}] \cdot
   \rho^*\alpha \cdot \phi^*a) \\
&\qquad + \psi_*(\rho^*\gamma \cdot \text{exp}(-[E_n])
   \cdot [X^{[n, n-1]}] \cdot  \rho^*\alpha \cdot \phi^*a) \\
&= \psi_*([X^{[n, n-1]}] \cdot \phi^*G(\gamma, n-1)
   \rho^*\alpha \cdot \phi^*a) \\
&\qquad + \psi_*((\text{exp}(-[E_n])
   \cdot [X^{[n, n-1]}]) \cdot \rho^*(\gamma\alpha) \cdot \phi^*a). \\
\endalign$$
Thus, $[\frak G(\gamma), \frak q_1(\alpha)](a)
= \psi_*((\text{exp}(-[E_n])
\cdot [X^{[n, n-1]}]) \cdot \rho^*(\gamma\alpha) \cdot \phi^*a)$
which is equal to $\text{exp}(\ad(\frak d))(\frak q_1(\gamma \alpha))(a)$
by the Lemma 3.9 in \cite{Leh}.
\endproof

\proposition{5.11} Assume that $n \in \Zee$ and $\gamma, \alpha \in H^*(X)$.
Then, the commutator $[\frak G(\gamma), \frak q_n(\alpha)]$ satisfies
the transfer property, i.e., we have
$$[\frak G(\gamma), \frak q_n(\alpha)]
= [\frak G(1_X), \frak q_n(\gamma \alpha)]
= [\frak G(\gamma \alpha), \frak q_n(1_X)]. \tag 5.12$$
\endstatement
\proof
Since the operator $\frak G(\gamma)$ is self-adjoint,
we see from (2.8) and (2.12) that we need only to prove (5.12)
for positive integers $n$. In the following,
we shall use induction on these positive integers $n$.
By Lemma 5.8, (5.12) is true for $n=1$. Assume that (5.12) is true for $n$.
Since $[\frak G(\gamma), \frak q_1(1_X)]
= \text{\rm exp}(\ad(\frak d))(\frak q_1(\gamma))$,
the same proof to (3.14) (replacing the operator $\frak f$ there
by $\frak G(\gamma)$) yields
$$[\frak G(\gamma), \frak q_{n+1}(\alpha)]
= - {1 \over n} \left \{ \sum_{j=0}^{+\infty} {1 \over j!} \cdot
    [\frak q_1^{(j+1)}(\gamma), \frak q_{n}(\alpha)]
- [[\frak G(\gamma), \frak q_{n}(\alpha)], \frak q_1'(1_X)] \right \}.$$
Now the transfer property of $[\frak G(\gamma), \frak q_{n+1}(\alpha)]$
follows from the induction hypothesis and the transfer property of
$[\frak q_1^{(j+1)}(\gamma), \frak q_{n}(\alpha)]$ in Proposition 3.17.
\endproof

\theorem{5.13} Let $k \in \Zee$ be nonnegative,
$\gamma \in H^*(X)$ be a nonzero cohomology class,
and $\epsilon(\frak G_{k}(\gamma)) = \frak G_{k}(\gamma)
+ W_0^{k+2}(\gamma)$. Then,
\par
{\rm (i)} for all $n_1, \ldots, n_{k+2} \in \Zee$ and
all $\alpha_1, \ldots, \alpha_{k+2} \in H^*(X)$,
$$[ \ldots [[\epsilon(\frak G_{k}(\gamma)), \frak q_{n_1}(\alpha_1)],
\ldots ], \frak q_{n_{k+2}}(\alpha_{k+2})] = 0;$$
\par
{\rm (ii)} for all $n_1, \ldots, n_{k+1} \in \Zee$ with
$\sum_{j = 1}^{k+1} n_j \ne 0$ and
all $\alpha_1, \ldots, \alpha_{k+1} \in H^*(X)$,
$$[ \ldots [[\epsilon(\frak G_{k}(\gamma)), \frak q_{n_1}(\alpha_1)],
\ldots ], \frak q_{n_{k+1}}(\alpha_{k+1})] = 0;$$
\par
{\rm (iii)} the leading term of $\frak G_{k}(\gamma)$
is $-W_0^{k+2}(\gamma)$;
\par
{\rm (iv)} $\frak G_{0}(\gamma) = -W_0^2(\gamma)$;
\par
{\rm (v)} $\frak G_{1}(\gamma) = -W_0^{3}(\gamma)$
if $K_X$ is numerically trivial. In particular,
$$\frak G_1(1_X) = \frak d = -W^3_0(1_X).  \tag 5.14$$
\endstatement
\proof
The first three statements follow from (5.3), (5.9) and Theorem 4.12.

To prove (iv), recall from the definitions that
$W_0^2(\gamma) = \frak L_0(\gamma)$.
By (5.6), we may assume that $\gamma$ is not a scalar multiple
of $1_X$. So $|\gamma| > 0$. Now,
$$[ \frak G_{0}(\gamma), \frak q_n(\alpha)] =
[ \frak G_0(1_X), \frak q_n(\gamma \alpha)]
= [ -\frak L_0(1_X), \frak q_n(\gamma \alpha)]  \tag 5.15$$
by (5.12) and (5.6). In view of Theorem 2.16 (ii), we have
$$[ \frak G_{0}(\gamma) +W^2_0(\gamma),\frak q_n(\alpha)] =
[-\frak L_0(1_X) ,\frak q_n(\gamma \alpha)]
+ [\frak L_0(\gamma) ,\frak q_n(\alpha)] = 0$$
for all $n \in \Zee$ and $\alpha \in H^*(X)$. Since $\Bbb H$ is irreducible,
$\frak G_{0}(\gamma) +W^2_0(\gamma)$ must be a scalar multiple of
the identity operator, which has to be zero since the bi-degree of
$(\frak G_{0}(\gamma) +W^2_0(\gamma))$ is $(0, |\gamma|)$.
So we have $\frak G_{0}(\gamma) = -W^2_0(\gamma)$.

Finally, to prove (v), we notice from Lemma 4.6 that
$$[ W^3_0(\gamma), \frak q_n(\alpha)] = -n W^2_n(\gamma \alpha)
= -n \frak L_n(\gamma \alpha).  \tag 5.16$$
Next, we have $\frak G_{1}(1_X) = \frak d$ in view of (5.7).
By (5.12) and Theorem 2.16 (iv),
$$[ \frak G_{1}(\gamma), \frak q_n(\alpha)] =
[ \frak G_1(1_X), \frak q_n(\gamma \alpha)]
= \frak q_n'(\gamma \alpha) = n \frak L_n(\gamma \alpha).  \tag 5.17$$
Combining this with (5.16), we see that
$[ \frak G_{1}(\gamma) +W^3_0(\gamma),\frak q_n(\alpha)] =0$
for all $n \in \Zee$ and $\alpha \in H^*(X)$.
As argued in the proof of (iv), we have
$\frak G_{1}(\gamma) = -W^3_0(\gamma)$.
\endproof

We remark that a result parallel to Theorem 5.13 (v) has been proved
by Frenkel and Wang \cite{F-W} within the framework of wreath products
(compare p.205 of \cite{Leh}). Also, $-W^{k}_0(\gamma)$ is
the degree-$0$ component of the vertex operator in (4.5).

\medskip\noindent
{\bf Definition 5.18.} Fix a positive integer $n$.
Define $\Bbb H_n'$ be the subring of $\Bbb H_n = H^*(\hil{n})$
generated by the following cohomology classes:
$$G_{i}(\gamma, n) = \frak G_{i}(\gamma)(1_{\hil{n}})
\tag 5.19$$
where $0 \le i < n$ and $\gamma$ runs over a linear basis of $H^*(X)$.

Note that the subring $\Bbb H_n'$ is generated by
$(n \cdot \dim H^*(X))$ elements.
Our goal is to show that $\Bbb H_n' = \Bbb H_n$, i.e.,
the cohomology ring $\Bbb H_n = H^*(\hil{n})$ is generated by
those $(n \cdot \dim H^*(X))$ classes $G_i(\gamma, n)$ in (5.19).
We shall use induction on the reverse lexicographic order $\prec$
of all the partitions $\mu = (\mu_1, \mu_2, \ldots)$ of $n$,
where $\mu_i$ denotes the number of parts in the partition $\mu$
equal to $i$ (compare with the proof of Theorem 4.10 in \cite{Leh}).
Under this ordering, the partition $(n, 0, ...)$ is the smallest.
Our induction goes as follows. First of all, we prove that
$$\prod_{j=1}^{n} \frak q_1(\alpha_{j}) |0\rangle
\in \Bbb H_n'  \tag 5.20$$
for all $\alpha_{j} \in H^*(X)$ in Lemma 5.23 below. Then, by assuming that
$$\prod_{i \ge 1} \prod_{j=1}^{\mu_i'} \frak q_i(\beta_{i,j})
|0\rangle \in \Bbb H_n'  \tag 5.21$$
for all $\beta_{i,j} \in H^*(X)$ and all $\mu' = (\mu_1', \mu_2', \ldots)
\prec \mu = (\mu_1, \mu_2, \ldots)$, we prove that
$$\prod_{i \ge 1} \prod_{j=1}^{\mu_i} \frak q_i(\alpha_{i,j})
|0\rangle \in \Bbb H_n'  \tag 5.22$$
for all $\alpha_{i,j} \in H^*(X)$.
Now we begin with the statement and proof of Lemma 5.23.

\lemma{5.23} For all cohomology classes $\alpha_{j} \in H^*(X)$, we have
$$\prod_{j=1}^{n} \frak q_1(\alpha_{j}) |0\rangle
\in \Bbb H_n'.  \tag 5.24$$
\endstatement
\proof
We shall use induction on $k$ to show that for all
$\alpha_{1}, \ldots, \alpha_k  \in H^*(X)$,
$$\prod_{j=1}^{k} \frak q_1(\alpha_{j}) \frak q_1(1_X)^{n-k}
|0\rangle \in \Bbb H_n'.  \tag 5.25$$

First of all, we claim that (5.25) is true for $k = 0$.
Indeed, in view of (5.6),
$$\frak q_1(1_X)^{n} |0\rangle
= n! \cdot 1_{\hil{n}} = (n-1)! \cdot G_0(1_X, n) \in \Bbb H_n'.$$

Next, assuming that (5.25) is true for some $k$ with $0 \le k < n$,
we shall verify that (5.25) holds as well if $k$ is replaced by $(k+1)$.
We may assume that $\alpha_{k+1}$ is homogeneous and $|\alpha_{k+1}|=s$.
By (5.9), we have $[\frak G_0(\alpha_{k+1}), \frak q_1(\alpha)] =
\frak q_1(\alpha_{k+1}\alpha)$. Taking the cup product of (5.25) with
$G_0(\alpha_{k+1}, n) \in \Bbb H_n'$, we obtain
$$\align
 \Bbb H_n'
&\ni G_0(\alpha_{k+1}, n) \cdot \left (
 \prod_{j=1}^{k} \frak q_1(\alpha_{j}) \frak q_1(1_X)^{n-k}
 |0\rangle  \right )    \\
&= \frak G_0(\alpha_{k+1})
 \prod_{j=1}^{k} \frak q_1(\alpha_{j}) \frak q_1(1_X)^{n-k}
 |0\rangle     \\
&= \sum_{i=1}^k \pm
 \frak q_1(\alpha_{1}) \cdots \frak q_1(\alpha_{i-1})
 \frak q_1(\alpha_{k+1} \alpha_{i}) \frak q_1(\alpha_{i+1})
 \cdots \frak q_1(\alpha_{k}) \frak q_1(1_X)^{n-k} |0\rangle \\
&\quad + (-1)^{s \sumlim{j=1}{k} |\alpha_j|} (n-k) \cdot
 \prod_{j=1}^{k} \frak q_1(\alpha_{j})
 \frak q_1(\alpha_{k+1}) \frak q_1(1_X)^{n-(k+1)}
 |0\rangle    \\
&\equiv (-1)^{s \sumlim{j=1}{k} |\alpha_j|} (n-k)
 \cdot \prod_{j=1}^{k+1} \frak q_1(\alpha_{j})
 \frak q_1(1_X)^{n-(k+1)} |0\rangle  \pmod {\Bbb H_n'}  \\
\endalign$$
where we have used the induction hypothesis in the first and last steps.
So (5.25) holds if $k$ is replaced by $(k+1)$.
This completes the proof of (5.24).
\endproof

\lemma{5.26} Fix $a, b$ with $1 \le a \le b$.
Let $\frak g \in \End(\Bbb H)$ be of bi-degree $(\ell, s)$, and
$$A = \frak q_{m_1}(\beta_{1}) \cdots
\frak q_{m_b}(\beta_{b}) |0\rangle.$$
Then, $\frak g(A)$ is equal to the sum of the following two terms:
$$\sum_{i=0}^{a-1} \sum_{\sigma_i} \pm
 \prod_{\ell \in \sigma_i^0} \frak q_{m_\ell}(\beta_{\ell})
 [\cdots [\frak g, \frak q_{m_{\sigma_i(1)}}(\beta_{\sigma_i(1)})],
 \cdots], \frak q_{m_{\sigma_i(i)}}(\beta_{\sigma_i(i)})]
 |0\rangle  \tag 5.27$$
and
$$\sum_{\sigma_a} (-1)^{\sum_{k=0}^{a-1}
 (s+\sum_{\ell=1}^k |\beta_{\sigma_a(\ell)}|)
 \sum_{\sigma_a(k)<j<\sigma_a(k+1)}|\beta_j|}
 \prod_{\ell \in \sigma_a^1} \frak q_{m_\ell}(\beta_{\ell}) \cdot$$
$$\cdot  [\cdots [\frak g,
 \frak q_{m_{\sigma_a(1)}}(\beta_{\sigma_a(1)})],
 \cdots], \frak q_{m_{\sigma_a(a)}}(\beta_{\sigma_a(a)})]
 \prod_{\ell \in \sigma_a^2} \frak q_{m_\ell}(\beta_{\ell})
 |0\rangle  \tag 5.28$$
where for each fixed $i$ with $0 \le i \le a$,
$\sigma_i$ runs over all the maps
$$\{\, 1, \ldots, i \,\} \to \{\, 1, \ldots, b \,\}$$
satisfying $\sigma_i(1) < \cdots < \sigma_i(i)$. Moreover,
$\sigma_i^0 = \{ \ell \, | \, 1 \le \ell \le b, \ell \ne
\sigma_i(1), \ldots, \sigma_i(i) \}$,
$\sigma_a^1 = \{ \ell \, | \, 1 \le \ell < \sigma_a(a),
\ell \ne \sigma_a(1), \ldots, \sigma_a(a) \}$,
and $\sigma_a^2 = \{ \ell \, | \, \sigma_a(a) < \ell \le b \}$.
\endstatement
\proof
Note that for all $i$ with $0 \le i < a$ and for all the above $\sigma_i$,
we move
$$[\cdots [\frak g, \frak q_{m_{\sigma_i(1)}}(\beta_{\sigma_i(1)})],
 \cdots], \frak q_{m_{\sigma_i(i)}}(\beta_{\sigma_i(i)})]$$
all the way to the right. This produces (5.27). In doing so,
we obtain (5.28) by repeatedly applying the elementary fact that
$$\frak g_1 \frak g_2 = [\frak g_1, \frak g_2] +
(-1)^{s_1 s_2} \frak g_2 \frak g_1  \tag 5.29$$
for two operators $\frak g_1, \frak g_2 \in \End(\Bbb H)$
of bi-degrees $(\ell_1, s_1)$, $(\ell_2, s_2)$ respectively.
\endproof

\theorem{5.30} For $n \ge 1$, the cohomology ring
$\Bbb H_n = H^*(\hil{n})$ is generated by
$$G_{i}(\gamma, n) = \frak G_{i}(\gamma)(1_{\hil{n}})$$
where $0 \le i < n$ and $\gamma$ runs over a linear basis of $H^*(X)$.
Moreover, the relations among these generators
are precisely the relations among the restrictions
$\frak G_i(\gamma)|_{\Bbb H_n}$ of the corresponding operators
$\frak G_i(\gamma)$ to ${\Bbb H_n}$.
\endstatement
\proof
Note that the second statement follows from the fact that
the operators $\frak G_i(\gamma)|_{\Bbb H_n}$ are defined
in terms of the cup products by the cohomology classes
$G_i(\gamma, n)$. In the following, we prove the first statement.

Let $\Bbb H_n'$ be defined as in Definition 5.18. We want to show that
$\Bbb H_n' = \Bbb H_n$. As indicated in the paragraph following
Definition 5.18, we use induction on the reverse lexicographic order
$\prec$ of all the partitions $\mu = (\mu_1, \mu_2, \ldots)$ of $n$.
By Lemma 5.23, (5.20) is true. In the following,
assuming that (5.21) is true for
all $\beta_{i,j} \in H^*(X)$ and all $\mu' = (\mu_1', \mu_2', \ldots)
\prec \mu = (\mu_1, \mu_2, \ldots)$, we shall prove (5.22).

Since $\mu = (\mu_1, \mu_2, \ldots) \ne (n, 0, \ldots)$,
we let $a$ be the smallest index such that $a > 1$ and $\mu_a \ge 1$
(i.e., $\mu = (\mu_1, 0, \ldots, 0, \mu_a, \mu_{a+1}, \ldots)$).
Also, some of the classes $\alpha_{1,j}$ in (5.22)
might be (scalar multiples of) $1_X$. Without loss of generality,
we may assume that $\alpha_{1,1} = \ldots = \alpha_{1, r-a} = 1_X$
for some $r \ge a$ and that $|\alpha_{1,j}| > 0$ for
all $j$ satisfying $(r-a) < j \le \mu_1$. Then (5.22) can be rewritten as
$$\frak q_1(1_X)^{r-a} \,
 \prod_{j=(r-a+1)}^{\mu_1} \frak q_1(\alpha_{1,j})
 \frak q_a(\alpha_{a,1}) \,
 \prod_{j=2}^{\mu_a} \frak q_a(\alpha_{a,j}) \,
 \prod_{i > a} \prod_{j=1}^{\mu_i} \frak q_i(\alpha_{i,j})
 |0\rangle \in \Bbb H_n'.   \tag 5.31$$
Put $\mu' = (a+\mu_1, 0, \ldots, 0, \mu_a-1, \mu_{a+1}, \ldots)$ and
$$\align
 A \,
&{\overset \text{def} \to =} \, \frak q_1(1_X)^r \,
 \prod_{j=(r-a+1)}^{\mu_1} \frak q_1(\alpha_{1,j}) \,
 \prod_{j=2}^{\mu_a} \frak q_a(\alpha_{a,j}) \,
 \prod_{i > a} \prod_{j=1}^{\mu_i} \frak q_i(\alpha_{i,j})
 |0\rangle\\
&{\overset \text{def} \to =} \,
 \frak q_{m_1}(\beta_{1}) \cdots \frak q_{m_b}(\beta_{b})
 |0\rangle    \tag 5.32
\endalign$$
where $b = (a+\mu_1)+(\mu_a-1)+ \mu_{a+1} + \ldots$ is
the length of $\mu'$. We have
$$\align
&m_1 = \ldots = m_{a+\mu_1}=1, m_i > 1
 \text{ for } (a+\mu_1) < i \le b,  \tag 5.33   \\
&\beta_1 = \ldots = \beta_{r}=1_X, |\beta_i| > 0
 \text{ for } r < i \le (a+\mu_1).  \tag 5.34   \\
\endalign$$
Since $A$ corresponds to the partition $\mu'$ and
$\mu' \prec \mu$, $A \in \Bbb H_n'$ by induction.
\medskip\noindent
{\bf Claim.} (5.31) is true as long as $|\alpha_{a,1}| = 4$.
\proof
Note that $a \le n$. So we see from the definition of
$\Bbb H_n'$ that $G_{a-1}(\alpha_{a,1}, n) \in \Bbb H_n'$.
Taking the cup product of $A \in \Bbb H_n'$ with
$G_{a-1}(\alpha_{a,1}, n) \in \Bbb H_n'$ yields
$$\Bbb H_n' \ni G_{a-1}(\alpha_{a,1}, n) \cdot A
= \frak G_{a-1}(\alpha_{a,1})(A).$$
Put $\frak e = \epsilon(\frak G_{a-1}(\alpha_{a,1}))
=\frak G_{a-1}(\alpha_{a,1})+ W_0^{a+1}(\alpha_{a,1})$. Then,
$$\Bbb H_n' \ni \frak G_{a-1}(\alpha_{a,1})(A)
= -W_0^{a+1}(\alpha_{a,1})(A) + \frak e(A).   \tag 5.35$$

Applying Lemma 5.26 to the operator $\frak e$,
we see that $\frak e(A)$ consists of two parts (5.27) and (5.28).
By (5.33), the number of $\frak q_1$'s in every nonvanishing term of
(5.27) is at least $(a+\mu_1) - (a-1) = \mu_1 + 1 > \mu_1$.
So every nonvanishing term in (5.27) corresponds to
some partition $\mu' \prec \mu$. By induction hypothesis,
(5.27) is contained in $\Bbb H_n'$.
By Theorem 5.13 (ii) (replacing the integer $k$ there by $(a-1)$),
we see that (5.28) is $0$. In summary,
$\frak e(A) \in \Bbb H_n'$. By (5.35), we obtain
$$-W_0^{a+1}(\alpha_{a,1})(A) \in \Bbb H_n'.   \tag 5.36$$

Now we apply Lemma 5.26 to the operator $-W_0^{a+1}(\alpha_{a,1})$.
So $-W_0^{a+1}(\alpha_{a,1})(A)$ consists of two parts (5.27) and (5.28).
Again, (5.27) is contained in $\Bbb H_n'$.
Let $N(\sigma_a)$ be the number of $\frak q_1$'s in
a nonvanishing term in (5.28) corresponding to $\sigma_a$.
By (5.33), $N(\sigma_a) \ge (a+\mu_1) - a = \mu_1$.
If $N(\sigma_a) > \mu_1$, then by induction hypothesis,
this nonvanishing term in (5.28) corresponding to $\sigma_a$
is contained in $\Bbb H_n'$. Also, $N(\sigma_a) = \mu_1$ if and only if
$m_{\sigma_a(1)} = \ldots = m_{\sigma_a(a)} =1$, i.e.,
$$1 \le \sigma_a(1) < \ldots < \sigma_a(a) \le (a+\mu_1)$$
by (5.33). So this nonvanishing term in (5.28) is of the form:
$$(-1)^{\sum_{k=0}^{a-1}
 (4+\sum_{\ell=1}^k |\beta_{\sigma_a(\ell)}|)
 \sum_{\sigma_a(k)<j<\sigma_a(k+1)}|\beta_j|}
 \cdot \prod_{\ell \in \sigma_a^1} \frak q_1(\beta_{\ell}) \cdot$$
$$\cdot  [\cdots [-W_0^{a+1}(\alpha_{a,1}),
 \frak q_1(\beta_{\sigma_a(1)})],
 \cdots], \frak q_1(\beta_{\sigma_a(a)})] \,
 \prod_{\ell \in \sigma_a^2} \frak q_{m_\ell}(\beta_{\ell})
 |0\rangle$$
which can be simplified to the following by Proposition 4.10 (i):
$$(-1)^{\sum_{k=0}^{a-1}
 (4+\sum_{\ell=1}^k |\beta_{\sigma_a(\ell)}|)
 \sum_{\sigma_a(k)<j<\sigma_a(k+1)}|\beta_j|}
 \cdot \prod_{\ell \in \sigma_a^1} \frak q_1(\beta_{\ell}) \cdot$$
$$\cdot (-1)^{a+1} \cdot
 \frak q_a(\alpha_{a,1}\beta_{\sigma_a(1)} \cdots \beta_{\sigma_a(a)})
 \, \prod_{\ell \in \sigma_a^2} \frak q_{m_\ell}(\beta_{\ell})
 |0\rangle.  \tag 5.37$$

Since $\alpha_{a,1} \in H^4(X)$, the term (5.37) being nonzero forces
$$|\beta_{\sigma_a(1)}|= \ldots =|\beta_{\sigma_a(a)}|=0. \tag 5.38$$
Since $1 \le \sigma_a(1) < \ldots < \sigma_a(a) \le (a+\mu_1)$,
we see from (5.34) that $1 \le \sigma_a(1) < \ldots < \sigma_a(a) \le r$.
So (5.37) can be further simplified to
$$\align
&\quad (-1)^{a+1} \cdot \prod_{\ell \in \sigma_a^1} \frak q_1(1_X)
 \frak q_a(\alpha_{a,1}) \, \prod_{\ell > \sigma_a(a)}
 \frak q_{m_\ell}(\beta_{\ell}) |0\rangle   \\
&= (-1)^{a+1} \cdot \frak q_1(1_X)^{|\sigma_a^1|}
 \frak q_a(\alpha_{a,1}) \frak q_1(1_X)^{r-a-|\sigma_a^1|} \cdot \\
&\qquad \cdot \prod_{j=(r-a+1)}^{\mu_1} \frak q_1(\alpha_{1,j}) \,
 \prod_{j=2}^{\mu_a} \frak q_a(\alpha_{a,j}) \,
 \prod_{i > a} \prod_{j=1}^{\mu_i} \frak q_i(\alpha_{i,j})
 |0\rangle     \\
&= (-1)^{a+1} \cdot \frak q_1(1_X)^{r-a} \,
 \prod_{j=(r-a+1)}^{\mu_1} \frak q_1(\alpha_{1,j}) \cdot \\
&\qquad \cdot \frak q_a(\alpha_{a,1}) \,
 \prod_{j=2}^{\mu_a} \frak q_a(\alpha_{a,j}) \,
 \prod_{i > a} \prod_{j=1}^{\mu_i} \frak q_i(\alpha_{i,j})
 |0\rangle.  \tag 5.39   \\
\endalign$$
Now there are exactly ${r \choose a}$ such terms in (5.28).
Therefore, by (5.36), we have
$${r \choose a} \cdot
 (-1)^{a+1} \cdot \prod_{j=1}^{\mu_1} \frak q_1(\alpha_{1,j})
 \frak q_a(\alpha_{a,1}) \,
 \prod_{j=2}^{\mu_a} \frak q_a(\alpha_{a,j}) \,
 \prod_{i > a} \prod_{j=1}^{\mu_i} \frak q_i(\alpha_{i,j})
 |0\rangle \in \Bbb H_n'.$$
It follows that (5.31) is true as long as $|\alpha_{a,1}| = 4$.
\endproof

We continue the proof of Theorem 5.30. In view of the above Claim,
it remains to prove that if there is an integer $s$
such that $0 \le s \le 3$ and (5.31) is true as long as
$|\alpha_{a,1}| \ge (s+1)$, then (5.31) holds as well
if $|\alpha_{a,1}| = s$. In view of Theorem 5.13, we put
$\w \frak e = \epsilon(\frak G_{a-1}(\alpha_{a,1}))
=\frak G_{a-1}(\alpha_{a,1})+ W_0^{a+1}(\alpha_{a,1})$. Then,
$$\Bbb H_n' \ni \frak G_{a-1}(\alpha_{a,1})(A)
= -W_0^{a+1}(\alpha_{a,1})(A) + \w \frak e(A).   \tag 5.40$$
As in the proof of the above Claim, we see from Lemma 5.26 that
$\Bbb H_n'$ contains
$$\sum_{\sigma_a} (-1)^{\sum_{k=0}^{a-1}
 (s+\sum_{\ell=1}^k |\beta_{\sigma_a(\ell)}|)
 \sum_{\sigma_a(k)<j<\sigma_a(k+1)}|\beta_j|}
 \cdot \prod_{\ell \in \sigma_a^1} \frak q_1(\beta_{\ell}) \cdot$$
$$\cdot (-1)^{a+1} \cdot
 \frak q_a(\alpha_{a,1}\beta_{\sigma_a(1)} \cdots \beta_{\sigma_a(a)})
 \, \prod_{\ell \in \sigma_a^2} \frak q_{m_\ell}(\beta_{\ell})
 |0\rangle   \tag 5.41$$
where $\sigma_a$ runs over all the maps
$\{\, 1, \ldots, a \,\} \to \{\, 1, \ldots, b \,\}$
with $\sigma_a(1) < \cdots < \sigma_a(a)$, and satisfies
$m_{\sigma_a(1)} = \ldots = m_{\sigma_a(a)} =1$. So by (5.33),
we have $\sigma_a(1) < \cdots < \sigma_a(a) \le (a + \mu_1)$.
Note that every nonvanishing term in (5.41) corresponds to
the partition $\mu$. Moreover,
if $|\beta_{\sigma_a(1)}|+ \cdots +|\beta_{\sigma_a(a)}| > 0$, then
$|\alpha_{a,1}|+|\beta_{\sigma_a(1)}|+ \cdots +|\beta_{\sigma_a(a)}|
\ge (s+1)$; so this nonvanishing term in (5.41) is already contained
in $\Bbb H_n'$ by our assumption. Therefore,
the subring $\Bbb H_n'$ contains
$$\sum_{\sigma_a} (-1)^{\sum_{k=0}^{a-1}
 (s+\sum_{\ell=1}^k |\beta_{\sigma_a(\ell)}|)
 \sum_{\sigma_a(k)<j<\sigma_a(k+1)}|\beta_j|}
 \cdot \prod_{\ell \in \sigma_a^1} \frak q_1(\beta_{\ell}) \cdot$$
$$\cdot (-1)^{a+1} \cdot
 \frak q_a(\alpha_{a,1}\beta_{\sigma_a(1)} \cdots \beta_{\sigma_a(a)})
 \, \prod_{\ell \in \sigma_a^2} \frak q_{m_\ell}(\beta_{\ell})
 |0\rangle   \tag 5.42$$
where $\sigma_a$ runs over all the maps
$\{\, 1, \ldots, a \,\} \to \{\, 1, \ldots, b \,\}$ with
$\sigma_a(1) < \cdots < \sigma_a(a) \le (a + \mu_1)$, and satisfies
$|\beta_{\sigma_a(1)}|= \cdots =|\beta_{\sigma_a(a)}|=0$.
So we have $\sigma_a(1) < \cdots < \sigma_a(a) \le r$ by (5.34).
Now as in the last paragraph (starting from the line below (5.37))
in the proof of the above Claim, we obtain
$${r \choose a} \cdot
 (-1)^{a+1} \cdot \prod_{j=1}^{\mu_1} \frak q_1(\alpha_{1,j})
 \frak q_a(\alpha_{a,1}) \,
 \prod_{j=2}^{\mu_a} \frak q_a(\alpha_{a,j}) \,
 \prod_{i > a} \prod_{j=1}^{\mu_i} \frak q_i(\alpha_{i,j})
 |0\rangle \in \Bbb H_n'.$$
So (5.31) holds if $|\alpha_{a,1}| = s$. This completes the proof
of Theorem 5.30.
\endproof


\Refs

\widestnumber\key{MMM}

\ref \key Bea \by A. Beauville
\paper Sur la cohomologie de certains espaces de modules de fibr\'es
vectoriels \jour Geometry and Analysis (Bombay, 1992), 37-40,
Tata Inst. Fund. Res., Bambay, 1995
\endref

\ref \key Bor \by R. Borcherds
\paper Vertex algebras, Kac-Moody algebras, and the Monster
\jour Proc. Natl. Acad. Sci. USA \vol 83 \yr 1986 \pages 3068--3071
\endref

\ref \key Che \by J. Cheah
\paper Cellular decompositions for nested Hilbert schemes of points
\jour Pac. J. Math.
\vol 183 \pages 39-90  \yr 1998
\endref

\ref \key dCM \by M. de Cataldo, L. Migliorini
\paper The Chow groups and the motive of the Hilbert scheme of points
on a surface \jour Adv. in Math. \vol 151\pages 283-312\yr 2000
\endref

\ref\key {EGL} \by G. Ellingsrud, L. G\" ottsche, M. Lehn
\paper On the cobordism class of the Hilbert scheme of a surface
\jour J. Alg. Geom. \vol 10 \yr 2001 \pages 81-100
\endref

\ref \key ES1 \by G. Ellingsrud, S.A. Str\o mme
\paper On the homology of the Hilbert scheme of points in the plane
\jour Invent. Math. \vol 87 \pages 343-352 \yr 1987
\endref

\ref \key ES2 \by G. Ellingsrud, S. A. Str\o mme
\paper Towards the Chow ring of the Hilbert scheme of $\Bbb P^2$
\jour J. reine angew. Math.  \vol 441 \pages 33-44  \yr 1993
\endref

\ref \key ES3 \by G. Ellingsrud, S. A. Str\o mme
\paper An intersection number for the punctual Hilbert scheme of a surface
\jour Transactions of A.M.S.  \vol 350 \pages 2547-2252  \yr 1999
\endref

\ref \key F-G \by B. Fantechi, L. G\"ottsche
\paper The cohomology ring of the Hilbert schemes of $3$ points on
a smooth projective variety
\jour J. reine angew. Math. \vol 439 \pages 147-158  \yr 1993
\endref

\ref \key Fog \by J. Fogarty
\paper Algebraic families on an algebraic surface
\jour Amer. J. Math.
\vol 90 \pages 511-520  \yr 1968 \endref

\ref \key FLM \by I. Frenkel, J. Lepowsky, A. Meurman
\paper Vertex operator algebra and the monster
\jour Pure and Applied Math.
\vol 134 \publ Academic Press  \yr 1988 \endref

\ref \key F-W \by I. Frenkel, W. Wang
\paper Virasoro algebra and wreath product convolution
\jour J. Alg. \vol 242\pages 656-671\yr 2001
\endref

\ref \key Go1 \by L. G\" ottsche
\book Hilbert schemes of zero-dimensional subschemes of smooth varieties
\bookinfo Lecture Notes in Mathematics {\bf 1572}
\publ Springer-Verlag \publaddr Berlin \yr 1994
\endref

\ref \key Go2 \by L. G\" ottsche
\paper On the motive of the Hilbert scheme of points on a surface
\jour Preprint
\endref

\ref \key Gro \by I. Grojnowski
\paper Instantons and affine algebras I:
the Hilbert scheme and vertex operators \jour Math. Res. Lett. \vol 3
\pages 275-291 \yr 1996
\endref

\ref \key Hai \by M. Haiman
\paper Hilbert schemes and Macdonald polynomials: the Macdonald
positivity conjecture \jour J. AMS. \vol 14\pages 941-1006\yr 2001
\endref

\ref \key Har \by  R. Hartshorne
\book Algebraic geometry
\publ Springer \publaddr Berlin-Heidelberg-New York\yr 1978
\endref

\ref \key Kac \by V. Kac
\book Vertex Algebras for Beginners, {\rm University Lecture Series}
{\bf 10} \publ AMS \publaddr Providence, Rhode Island \yr 1997
\endref

\ref \key K-T \by S.~Kumar, J.~F.~Thomsen
\paper Frobenius
splitting of Hilbert schemes of points on surfaces
\jour Preprint
\endref

\ref \key Leh \by M. Lehn
\paper Chern classes of tautological sheaves on Hilbert schemes of
points on surfaces \jour Invent. Math. \vol 136 \yr 1999
\pages 157-207 \endref

\ref\key LS1\by  M. Lehn and C. Sorger,
\paper Symmetric groups and
the cup product on the cohomology of Hilbert schemes
\jour Duke Math. J. \vol 110 \pages 345--357\yr 2001
\endref

\ref \key LS2 \by  M. Lehn and C. Sorger,
\paper The cup product of the Hilbert scheme for $K3$ surfaces
\jour Preprint, math.AG/0012166
\endref

\ref \key LQZ \by W.-P. Li, Z. Qin, Q. Zhang
\paper On the geometry of the Hilbert schemes of points
in the projective plane  \jour Preprint
\endref

\ref \key Mar \by E. Markman
\paper Generators of the cohomology ring of moduli spaces of sheaves
on symplectic surface \jour Preprint
\endref


\ref \key Na1 \by H. Nakajima
\paper Heisenberg algebra and Hilbert schemes of points on
projective surfaces \jour Ann. Math. \vol 145 \yr 1997 \pages 379-388
\endref

\ref \key Na2 \by H. Nakajima
\book Lectures on Hilbert schemes of points on surfaces,
{\rm University Lecture Series,} {\bf 18}
\publ American Mathematical Society \publaddr Providence, RI \yr 1999
\endref

\ref \key Tik \by A. S. Tikhomirov
\paper The variety of complete pairs of zero-dimensional subschemes of
an algebraic surface \jour Izv. Math.
\vol 61 \pages 1265-1291  \yr 1997  \endref

\ref \key V-W \by C. Vafa, E. Witten \paper A strong coupling test
of $S$-duality \jour Nucl. Phys. B \vol 431 \yr 1995 \pages 3-77
\endref

\ref \key Wan \by W. Wang
\paper Hilbert schemes, wreath products, and the McKay correspondence
\jour Preprint, math.AG/9912104
\endref

\ref \key Yos \by K. Yoshioka
\paper Betti numbers of moduli of stable sheaves on some surfaces.
$S$-duality and Mirror Symmetry (Trieste, 1995)
\jour Nuclear Phys. B Proc. Suppl.
\vol 46 \pages 263-268  \yr 1996  \endref


\endRefs
\enddocument